\documentclass[12pt]{article}
\usepackage{graphicx}
\usepackage[bookmarks=false]{hyperref}
\usepackage{amsmath}
\usepackage{amssymb}

\def\A{\mathbb{A}}
\def\R{{\mathbb R}}

\def\Z{{\mathbb Z}}
\def\N{{\mathbb N}}

\def\af{\alpha}

\def\ds{\displaystyle}

\def\all{\forall}
\def\lam{\lambda}

\def\sl2{SL(2,\R)}
\def\S{\mathbb{S}^1}
\def\Oph{\mbox{Homeo}^+(\S)}
\def\cOph{\widetilde{\mbox{Homeo}^+(\S)}}

\def\csl2{\widetilde{\sl2}}

\def\pd{$(+d)$}
\def\x{{\bf x}}
\def\y{{\bf y}}
\newcommand{\bea}{\begin{eqnarray}}
\newcommand{\eea}{\end{eqnarray}}
\newcommand{\beq}{\begin{equation}}
\newcommand{\eeq}{\end{equation}}
\newcommand{\barr}{\begin{array}}
\newcommand{\earr}{\end{array}}
\newcommand{\ovl}{\overline}
\newtheorem{theo}{Theorem}[section]
\newtheorem{prop}{Proposition}[section]
\newtheorem{defn}{Definition}[section]
\newtheorem{lemma}{Lemma}[section]
\newtheorem{cor}{Corollary} [section]

\newtheorem{rmk}{Remark}[section]
\def\bth{\begin{theo}}
\def\eth{\end{theo}}
\def\bpr{\begin{prop}}
\def\epr{\end{prop}}
\def\bdf{\begin{defn}}
\def\edf{\end{defn}}
\def\brmk{\begin{rmk}}
\def\ermk{\end{rmk}}

\begin{document}

\title{Twist number and order properties of periodic orbits}%
\author{Emilia Petrisor}
\date{}
\maketitle

\begin{abstract}
A less studied numerical
characteristic  of  periodic orbits of area preserving twist maps of the annulus is the twist  or torsion number, called initially the amount of rotation \cite{mather}. It  measures  the average
rotation of tangent vectors under the action of the derivative of the map
along that orbit, and characterizes  the degree of complexity of the dynamics.\par
The aim of this paper is to give new insights into the definition and  properties
of the twist  number,  and to relate  its
  range  to the order properties of periodic orbits. We derive an algorithm
to deduce the exact value or a demi--unit interval containing the exact value of the twist number.\par
 We prove that  at a period--doubling
bifurcation threshold of a mini-maximizing periodic orbit, the new born doubly periodic orbit has the absolute  twist number larger than
  the absolute twist  of the   original  orbit after bifurcation. We give examples of periodic orbits having large absolute twist number, that are badly ordered, and illustrate how characterization of these orbits only by their residue can lead to incorrect results.\par
  In connection to the study of the twist number of periodic orbits of standard--like maps we introduce a new tool, called $1$--cone function.
  We prove that the location of  minima of this function with respect to the vertical  symmetry lines of a   standard--like map   encodes a  valuable information on the symmetric periodic orbits and their twist number.
\end{abstract}

\section{Introduction}\label{introTwist}
The study of the dynamics of area preserving positive twist maps of the annulus is mainly concerned
with the  characterization of its invariant sets, and  the dynamical behaviour of a map restricted to  such sets.
Among invariant sets, rotational periodic orbits have been thoroughly studied and classified according to their linear stability, extremal type, order properties  (see \cite{meiss92} for a survey).  A numerical characteristic associated to a  periodic orbit is the rotation number, which measures the average rotation of the orbit around the annulus.  In \cite{mather} Mather defined also the amount of rotation, which is called twist number in \cite{angenent}
or torsion number in \cite{crovis}. The twist number  of a   periodic orbit characterizes  the average rotation of tangent vectors under the action of the tangent  map
along the map's orbit.  Mather  related the twist number with the Morse index of the corresponding critical sequence (configuration), and Angenent \cite{angenent} proved that in the space of $(p,q)$--sequences  a   critical point of the $W_{pq}$--action, corresponding to a periodic orbit of twist number greater than 0 is  connected by the negative gradient flow of the action, through a heteroclinic connection, with a
  sequence corresponding to an orbit of zero twist number.  Using topological arguments, Crovisier  \cite{crovis}  proved the existence of orbits of zero--torsion (twist) number for any  real number in the rotation number  set  of a twist map of the annulus (not necessarily area preserving).  \par

 In order to detect new classes of dynamical systems exhibiting periodic orbits of non--zero twist (torsion),   B\'{e}guin and Boubaker gave \cite{beguin} conditions ensuring that some  area preserving
diffeomorphisms of the disc $\mathbb{D}^2$, and particular diffeomorphisms of the torus $\mathbb{T}^2$  exhibit  such orbits.\par

In Mather's and Angenent's definitions of the twist number this a positive quantity. However its natural definition leads to a signed number
\cite{crovis}. The twist number of periodic orbits of positive twist maps is zero or negative, while for those of negative twist maps it is zero or positive. In our approach we keep its  natural sign and call absolute twist, the absolute   value of the twist number. \par

 During almost 30 years  since the definition of the twist number was given, no  periodic orbit of absolute twist number greater than $1/2$ was detected in the dynamics of classical twist maps (standard map, Fermi--Ulam map, etc). That is why we wonder under what conditions  can such orbits appear. \par
 From Aubry--Mather theory \cite{aubry}, \cite{matherIHES} it is known that non--degenerate minimizing and associated mini--maximizing $(p,q)$--type periodic orbits are well ordered and the minimizing orbit has zero twist number. Thus it is natural to investigate whether $(p,q)$--periodic orbits that are not given by this theory can be well--ordered  or not and how order properties are related to the value of their twist number. \par

 Our starting point was not  however  the study of the twist number. In an attempt to characterize  the dynamical behaviour of twist maps after the breakdown of the last KAM invariant torus (the case of standard map) or in a half annulus where no invariant circle exist (the case of Fermi--Ulam map \cite{korsch}, or tokamap \cite{balescu}),  we identified  in the phase space of such  maps, special regions where no minimizing periodic orbits can land. Instead we noticed that any periodic orbit having at least two points in such a region is badly ordered and typically has absolute twist greater than  $1/2$.\par

 Thus we were led to a
deeper study of the  twist number of periodic orbits, and their
properties. In this paper we complement the results on twist number
reported in \cite{mather}, \cite{angenent}, and show that periodic
orbits of non--zero twist numbers are typically born through a period
doubling bifurcation.  We give examples of periodic orbits  that
have large absolute twist number. As far as we know no such orbits  were identified  before in the study of the most known twist maps, standard--like maps. Some of these orbits  are ordered and other are unordered.
In order to explain why some twist maps can exhibit a sequence of bifurcations of a periodic orbit leading to  increasing  the absolute  twist number up to its maximal value
we introduce two subclasses of standard--like maps, USF maps and TSF maps. \par

The paper is organized as follows. In Section \ref{btwist} we give
the basic properties of twist area preserving diffeomorphisms of the
annulus, relevant to our study. In order to derive in Section \ref{twistnosect} sufficient conditions that favor  the existence of  ordered periodic orbits of large absolute twist number,  we set up in Section \ref{TwistFolding} a framework of study, defining the class of twist maps exhibiting the so called strong folding property, and a function, called 1--cone function.

 The $1$--cone function is defined on the phase--space of a twist map, and takes negative values within the region where the map exhibits strong folding property. We prove that the  restriction of this function to a periodic orbit gives information on the eigenvalues of the Hessian matrix
associated to that orbit (Lemma \ref{lemmaeigendelta}).  Analyzing the 1--cone function associated to a standard--like map, we show  that such maps can have either  a connected strong folding region including one of the vertical symmetry lines or a two--component strong folding region including both vertical symmetry lines. These results will be exploited in Section 5 to explain the variation of the twist number.\par

In Section \ref{twistnosect} we revisit the
definition and properties of the twist number of a periodic orbit
based on the structure of the universal covering group of the group
$\sl2$ (a system of coordinates on this group allowing to decipher
its topology is presented in \ref{apendix2}). The twist number is
defined as the translation number of a circle map induced by the
monodromy matrix associated to the periodic orbit.  We review in
\ref{apendix1} the properties of the translation number of an
orientation preserving homeomorphism of the unit circle and point
out the particularities of the translation number of circle
homeomorphisms induced by matrices in $\sl2$.\par
 Theorem 4.2 gives the relationship between the twist number value of a $(p,q)$--periodic orbit, and the position of the real number $0$ with respect to the sequence of interlaced eigenvalues of the Hessian matrix $H_{q}$, associated to the corresponding
  $(p,q)$--sequence, and of a symmetric matrix derived from $H_{q}$.
   This theorem complements results from \cite{mather} and \cite{angenent}.

   Based on this theorem we derive an algorithm to deduce the exact value of the twist number or a demi--unit interval that contains the twist number.

 The main result in Section \ref{twistnosect} is the
  Proposition \ref{propperiodD},  which shows that  periodic orbits of large absolute twist number are born through a period doubling bifurcation. More precisely it states that if
  at some
threshold,  a period doubling bifurcation of a
mini--maximizing $(p,q)$--periodic orbit occurs, with transition
from  elliptic to inverse hyperbolic orbit,  then the new born
$2q$--periodic orbit is elliptic, having the twist number within the
interval $(-1,-1/2)$.  \par

In Section \ref{TwistOrderP} we give examples of periodic orbits that  have large absolute twist number and are unordered.
  A natural question is whether a positive twist map can also exhibit   ordered $(p,q)$--periodic orbits ($p,q$ relative prime integers),  of twist number less than $-1/2$ (we note that periodic orbits of large absolute twist number, born through a period doubling bifurcation are of type $(2p, 2q)$).  We give a positive answer to this question, giving an example of  three--harmonic standard map
having an ordered  $(1,2)$--periodic orbit, of  twist number $\tau=-1$.   \par
In order to illustrate the different range of the twist number  of a periodic orbit in the three-harmonic standard map in comparison to that of  a periodic orbit of the same type of the standard map, we deduce (Proposition \ref{PropUSFbeh}) the bifurcations that a periodic orbit of a family of standard--like maps, with a uni-component strong folding region, undergo necessarily.

 The concrete examples of periodic orbits of large absolute twist number, given in this section,  illustrate that the classical method used for more than 30 years in classification of periodic orbits can lead to an erroneous conclusion.\par

\section{Background on twist maps}\label{btwist}
We recall basic properties of twist maps, relevant to our approach.
For a detailed presentation of classical results concerning dynamics
of this type of maps, the reader is referred to
\cite{gole2001}, \cite{meiss92}.
 \par Let $\mathbb{S}^1=\R/\Z$ be
the unit circle,   $\mathbb{A}=\mathbb{S}^1\times \R$, the infinite
annulus , and  $\pi:\R^2\to\mathbb{A}$, the covering projection,
$\pi(x,y)=( x\,\,\mbox{mod}\,\, 1, y)$.
\par

We consider a $C^1$--diffeomorphism $F:\R^2\to\R^2$,
$F(x,y)=(x',y')$, satisfying the following properties:\par i) $F$ is
exact area preserving map, isotopic to the identity;\par
 ii)
$F(x+1,y)=F(x,y)+(1,0)$, for all $(x,y)\in\R^2$;\par
 iii) $F$ has
uniform positive twist property, i.e. $\displaystyle\frac{\partial
F_1}{\partial y}\geq c>0$.\par

The map $F$ defines a $C^1$--diffeomorphism, $f:\A\to\A$, such that
$\pi\circ F=f\circ \pi$. Both maps $f$ and  $F$ are called  {\it area
preserving positive twist maps} or simply twist maps. $F$ is a lift of $f$.\par
  In the sequel we will
switch from $f$ to $F$, and conversely,  without comment.\par

 The motion of a point
around the annulus is characterized by its  rotation number. The
orbit of a point $z\in\A$ has a rotation number if there exists the
limit:
\beq\label{RotNo}\rho=\lim_{n\to\infty}\displaystyle\frac{x_n-x}{n},\eeq
where $(x,y)\in\R^2$ is a lift of $z$, and $(x_n,y_n)=F^n(x,y)$. The
rotation number  does not depend on   the chosen
point $z$ on the orbit or the lift $(x,y)$. For different lifts of the map $f$, the corresponding rotation numbers differ by an integer.\par
 Let $p,q$ be two relative prime
integers, $q>0$. A $(p,q)$--type orbit of the twist map, $F$, is an
orbit $(x_n,y_n)=F^{n}(x_0, y_0)$, $n\in\Z$, such that
$x_{n+q}=x_n+p$, for every $n\in\Z$.  The $\pi$--projection of such
an orbit onto the annulus is called  $(p,q)$-- periodic orbit or
simply $q$-periodic orbit of the map $f$, and its rotation number is  $p/q$.
\par The linear stability of a $(p,q)$--periodic orbit, $\mathcal{O}(z)$, is
characterized by the multipliers of the derivative, $D_zf^q\in
SL(2,\R)$, which in turn  are determined by the trace,
$\mbox{tr}(D_zf^q)$.   Instead of  trace one can also use the residue \cite{greene}
to classify the $q$--periodic orbits. Namely, the residue is defined
by: \beq\label{residueF} R=\displaystyle\frac{
2-\mbox{trace}(D_zf^q)}{4},\eeq where $q$ is the least period of the
orbit.\par

 The dynamics of a twist map, $F$, has a variational
formulation, that is the $F$--orbits are associated to critical
points of an action, and vice--versa. We discuss only the case of $(p,q)$--type orbits.\par
More precisely,  $F$ admits a  $C^2$--generating function  $h:\R^2\to\R$ (unique up to additive constants), such that
$h_{12}(x,x')=\partial h/\partial x\partial x'<0$ on $\R^2$, and
 $F(x,y)=(x', y')$ iff $y=-h_1(x,x')$, $y'=h_2(x,x')$  ($h_1=\partial h/\partial x, h_2=\partial h/\partial x')$.

The generating function defines an action
$W_{pq}$, on the space
 $\mathcal{X}_{pq}$, of
$(p,q)$--type sequences of real numbers, ${\bf x}=(x_n)_{n\in\Z}$,
i.e. sequences satisfying the property, $x_{n+q}=x_n+p$, for all
$n\in\Z$. The space  $\mathcal{X}_{pq}$  is identified with the affine
subspace of $\R^{q+1}=\{(x_0, x_1, \ldots, x_q)\}$, of equation
$x_q=x_0+p$, and the action is defined by: $$
W_{pq}({\bf x})=\sum_{k=0}^{q-1}h(x_k, x_{k+1})$$ There is a one to
one correspondence between critical points,  ${\bf x}=(x_n)$, of the $W_{pq}$--action,
and $(p,q)$--type orbits, $(x_n, y_n)$ of the positive twist map,
$F$:  \beq\label{criticseqpq} (x_n, y_n)=(x_n, -h_1(x_n, x_{n+1})),  \quad \all\,\, n\in\N\eeq

By Aubry--Mather theory \cite{aubry}, \cite{matherIHES}, for each  pair of relative prime integers
$(p,q)$, $q>0$, a twist area preserving map  has at least two periodic orbits  of
type $(p,q)$. One corresponds to a non--degenerate  minimizing
$(p,q)$--sequence of the action, and by abuse of language it is also
called in the following,
 minimizing $(p,q)$--type orbit.   Associated
to such an orbit there is a second one, corresponding to a
mini--maximizing sequence, and  called mini--maximizing
$(p,q)$--type periodic orbit.\par

Let $H_{q}$ be the  Hessian matrix of the action, $W_{pq}$, at a
critical point ${\bf x}=(x_n)$. The signature (the number of
negative and positive eigenvalues) of the Hessian $H_q$, at a
non--degenerate minimizing sequence is $(0, q)$, while at the
corresponding mini--maximizing sequence it is $( 1,q-1)$. The number
of negative eigenvalues is called the Morse index of the critical
sequence.
\par

A salient property of  a non--degenerate minimizing
$(p,q)$--periodic orbit of a positive twist map, $f$,  is that it is
well-ordered.
\par An invariant set $M$ of a positive twist map, $f$, is well
ordered if  for every $(x,y), (x', y')'\in\pi^{-1}(M)$, we have $x<x'$ iff
$F_1(x,y)<F_1(x',y')$, where $F_1$ is the first component of the
lift $F$.

A well--ordered $(p,q)$--orbit is also called Birkhoff orbit, while
a badly--ordered orbit is called non--Birkhoff.\par
 A
$(p,q)$--periodic orbit of the twist map $f$ is well ordered iff the
coresponding $(p,q)$--sequence $\x$ is cyclically ordered, i.e.
\beq\label{COorder} \x\leq \tau_{ij}\x\quad \mbox{or}\quad
\tau_{ij}\x\leq \x, \quad\all\,\, i,j\in\Z,\eeq
where $\tau_{ij}$ is the translation defined by:
$$(\tau_{ij}\x)_k=x_{k+i}+j, \quad \all\,\, \x=(x_k)\in\R^\Z, i,j,k\in\Z$$
By Aubry--Mather theory $W_{pq}$--minimizing sequences are
 cyclically ordered. Moreover, the mini--maximizing sequences are also
 ordered with respect to minimizing sequences \cite{matherIHES}. \par

The space $\mathcal{X}_{pq}$ of $(p,q)$--sequences is
partially ordered with respect to the order inherited from $\R^\Z$.
Namely if $\x=(x_k)$, $\y=(y_k)\in\R^\Z$, then:
\beq\label{seqorder}\barr{ll} \x \leq \y&\Leftrightarrow\quad
x_k\leq y_k, \all\,\, k\in\Z\earr\eeq

One also  defines the following relations: \beq\label{relless}\barr{ll}\x
<\y&\Leftrightarrow\quad \x \leq \y, \mbox{but}\,\, \x\neq
\y\\
\x\prec \y&\Leftrightarrow\quad  x_k<y_k, \all\,\, k\in\Z\earr \eeq

 To a $(p,q)$--sequence $\x=(x_k)$ one associates the Aubry diagram, i.e the graph
of the
piecewise affine function $A:\R\to\R$, that interpolates linearly
the points  $(k, x_k), k\in\Z$.  The Aubry function associated to a minimizing $(p,q)$--sequence is increasing. If the Aubry diagram of a sequence
$\y=(y_k)$ do not intersect the Aubry diagram of the minimizing sequence $\x=(x_k)$, then the two sequences are comparable, i.e. either $\x\prec \y$ or $\y\prec \x$. If the two diagrams intersect transversally, then the sequence $\y=(y_k)$ corresponds to an unordered $(p,q)$--periodic orbit.
\par
The most studied twist maps, both theoretically and numerically, are the standard--like maps.
 A standard--like map is a    twist map, $F_\epsilon$,   defined by a Lagrangian generating function of the form $h(x, x')=\dfrac{1}{2}(x-x')^2-\epsilon V(x)$, where  $V$ is a fixed 1--periodic even function of class $C^r$, $r\geq 3$, and $\epsilon\in\R$ is  a parameter:
$$\barr{lll} x'&=& x+y-\epsilon V'(x)\\
y'&=&y-\epsilon V'(x),\earr$$
Classical standard map (\ref{stdmapdef}) corresponds  to the potential $V(x)=-\ds\frac{1}{(2\pi)^2}\cos(2\pi x)$.\par
The twist map,  $F_\epsilon$, is   reversible, i.e. it factorizes as  $F_\epsilon= I\circ R$, where $R$ and $I$ are the  involutions, $R(x,y)=(-x, y-\epsilon V'(x))$, $I(x,y)=(-x+y, y)$. The symmetry lines, that are the fixed point sets $\mbox{Fix}(R)$, $\mbox{Fix}(I)$, have both two
 components: $\Gamma_0: x=0$,  $\Gamma'_0=0.5$,  respectively  $\Gamma_1: x=y/2$, $\Gamma'_1: x=(y+1)/2$.
 The $R$--invariant  orbits   are called symmetric orbits.\par

As we shall show in Section 4, the twist number of a $(p,q)$--periodic orbit one deduces analyzing the  relative position with respect to $0$ of the eigenvalues of the associated Hessian matrix and of a companion symmetric matrix.
 In the following section we prove that diagonal entries of these matrices encode information on the position of the smallest eigenvalue.
 In order to extract such an  information
  we introduce a new tool of study of a twist map, the 1--cone function, whose restriction to $(p,q)$--periodic orbit gives the diagonal entries in the two symmetric matrices.  Depending on the location of the minima of the $1$--cone function,   with respect to the vertical symmetry--lines of a  standard--like map we define two subclasses of this set of twist maps: the class of USF maps, respectively TSF maps. Then in Section 5 we prove that for a family of USF maps, respectively TSF maps,  one can predict  the sequence of bifurcations that  their  symmetric periodic orbits can undergo, as well as the variation of their twist number, only from the analysis of the diagonal entries of the Hessian and the companion matrix.

\section{$1$--cone function and  strong folding property of a twist map}\label{TwistFolding}

  Let $f$ be a positive or negative twist APM defined on the annulus $\mathbb{A}$.
The derivative of the map at an arbitrary point has the entries denoted as follows:
$$D_zf=\left[\barr{cc} a(z)&b(z)\\c(z)&d(z)\earr\right]$$
The  twist APM,  $f$, is called $(+d)$--map if $d(z)>0$, $\all
z\in\mathbb{A}$, respectively $(+a)$--map, if $a(z)>0$, $\all
z\in\mathbb{A}$. The two conditions are not equivalent.
\par
We study only twist maps $f$ that are $(+d)$-- maps, because  their inverses are  $(+a)$--maps. \par
  The \pd--condition ensures that the lines
of the vertical foliation of the annulus are mapped to graphs of
functions over the action, while the sign of $a(z)$
illustrates how the map acts on the horizontal foliation of the
annulus.   If for a \pd--twist map the function $a(z)$ is negative on a segment of an horizontal circle,
$\mathbb{S}^1\times \{y\}$, then the $f$--image  of that circle    is
not a graph over $x$. We say that the circle is folded by
the map.

\bdf A $(+d)$--twist map is said to  satisfy the folding property if
there exist regions in the phase space where $a\leq 0$. The subset
$\mathcal{F}=\{z\in\mathbb{A}\,|\, a(z)\leq 0\}$
is called  {\it folding region}.\edf
\par   The  standard map is a \pd--twist APM which
 exhibits folding property  for $\epsilon\geq 1$. The vertical
strip: \beq\label{foldstd} \mathcal{F}_\epsilon=\{(x,y)\,|\,
-\ds\frac{1}{2\pi}\arccos(1/\epsilon)\leq x\leq  \ds\frac{1}{2\pi}\arccos(1/\epsilon)\},
\quad \epsilon\geq 1\eeq is its folding region.

\par In the
sequel we refer only to \pd--positive twist maps, and all results can be
modified accordingly for \pd--negative twist ones.  In order to reveal the  properties of  twist maps exhibiting a folding region,
we use geometrical arguments, that are more intuitive and allow  visualizations   in the phase space. The same properties could be
presented in configuration space in a more formal setting.\par
We associate to each point of the phase space a pair of vectors:
\begin{equation}\begin{array}{l}\label{conesDef}
\barr{lll}v^+_1(z)&=&D_{f^{-1}(z)}f(e_2))=(b(f^{-1}(z)), d(f^{-1}(z)))^T,\\
v^-_1(z)&=&D_{f(z)}f^{-1}(-e_2)=(b(z), -a(z))^T,\earr
\end{array}\end{equation}
where $e_2=(0,1)^T$ is the vertical tangent vector at  $f^{-}(z)$, respectively at $f(z)$.
The
directions   to which the vectors   $v_1^+(z)$, $v_1^-(z)$  are pointing  illustrate how the
linear maps $D_{f^{-1}(z)}f$, $D_{f(z)}f^{-1}$,  deviate the the upward vertical
through $f^{-1}(z)$, respectively the
downward vertical through $f(z)$.  If at a point the slopes of the two vectors, $s_1^+(z)=\ds\frac{d}{b}(f^{-1}(z))$, $s_1^-(z)=-\ds\frac{a}{b}(z)$, are related by
$s_1^-(z)\geq s_1^+(z)$, one says that
  $1$--cone crossing occurs at $z$. It is known that  no minimizing orbit can pass through such a point    \cite{kayperciv}. \par
In a folding region, since $a(z)\leq 0$, the two slopes are both positive, and $1$--cone crossing is possible.
It is easy to check that 1--cone crossing can occur in the phase
space of a \pd--twist map only if the map has folding property, and
the set, $$\rm{Cone}_1=\{z\in \mathcal{\mathbb{A}}\,|\,
s_1^-(z)\geq s_1^+(z)\},$$  of the points where 1--cone crossing occurs
is necessarily included in the folding region.\par A\pd--twist map
with the property that the $1$--cone crossing occurs in a subset of a folding region   is said to exhibit  {\it strong folding property}, and the set $\rm{Cone}_1$ is called {\it strong folding region}.
\par

\par The  standard map   exhibits a  strong folding region
 for each $\epsilon\geq 2$:
\beq\label{1conestd}\rm{Cone}_1(\epsilon)=\{(x, y)\in
\mathcal{F}\,|\,-\ds\frac{1}{2\pi}\arccos{(2/\epsilon)}\leq x\leq
\ds\frac{1}{2\pi}\arccos{(2/\epsilon)}\}\eeq

\par

In order to show how the strong folding property influences the dynamical behaviour of a twist map, we  define a real function $\delta$, which associates to each point $z=(x,y)$ the difference between the slopes $s_1^+(z)$, $s_1^-(z)$:
 $$
 \delta(x,y)=\dfrac{d(x_{-1}, y_{-1})}{b(x_{-1}, y_{-1})}+\dfrac{a(x,y)}{b(x,y)},\quad (x_{-1},y_{-1})=f^{-1}(x,y)$$
We call this function, $1$--cone function associated to the twist map.
 $\delta$ takes negative values at the points of a strong folding region of a \pd--twist map.\par
\begin{figure}
\centerline{\scalebox{0.7}{\includegraphics{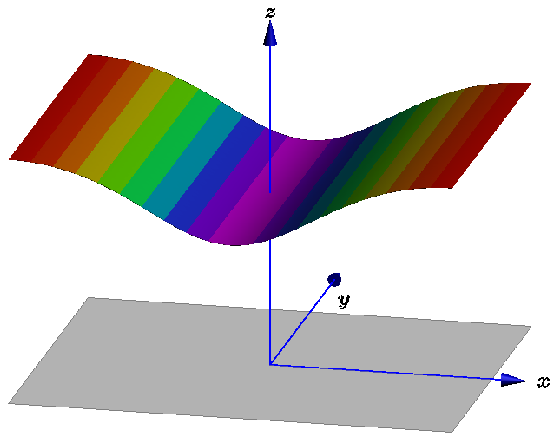}\hspace{0.25cm}\includegraphics{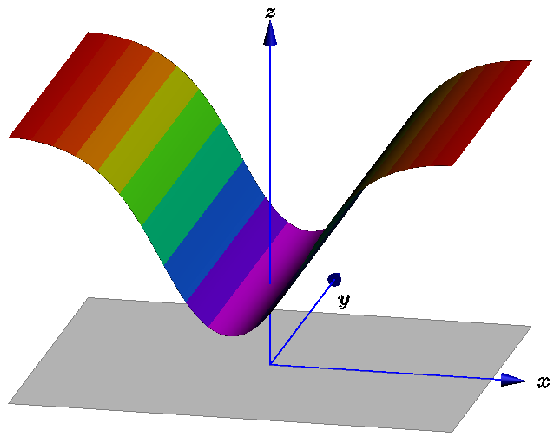}\hspace{0.25cm}\includegraphics{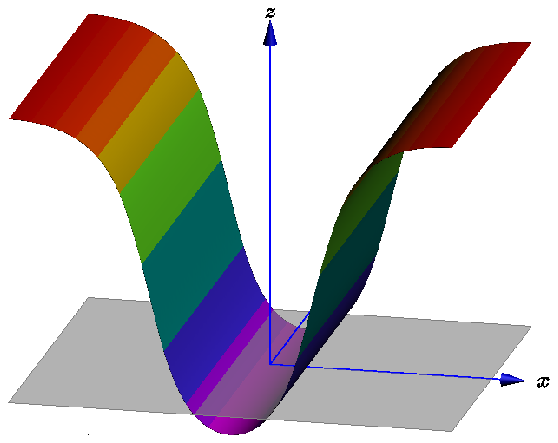}}}
\caption{\label{deltastd}  The graph of the function $\delta$ with respect to the plane $xOy$, associated to the standard map corresponding respectively to $\epsilon=0.5, 1.2, 2.3$.}
\end{figure}
\begin{figure}
\centerline{\scalebox{0.7}{\includegraphics{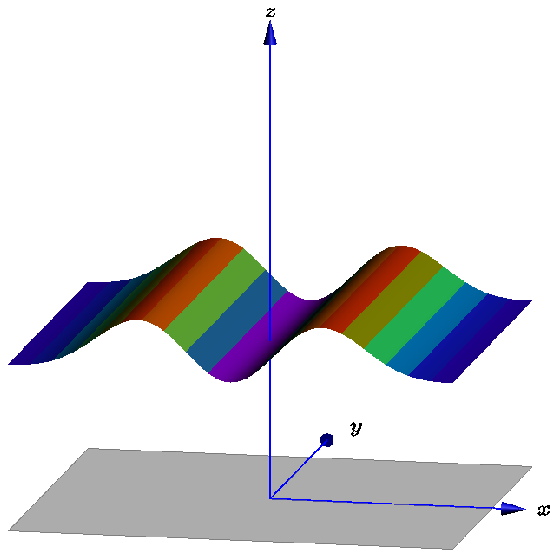}\hspace{0.25cm}\includegraphics{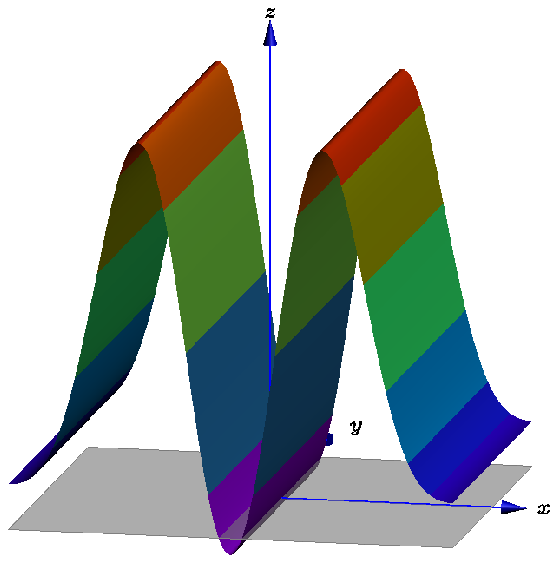}\hspace{0.25cm}\includegraphics{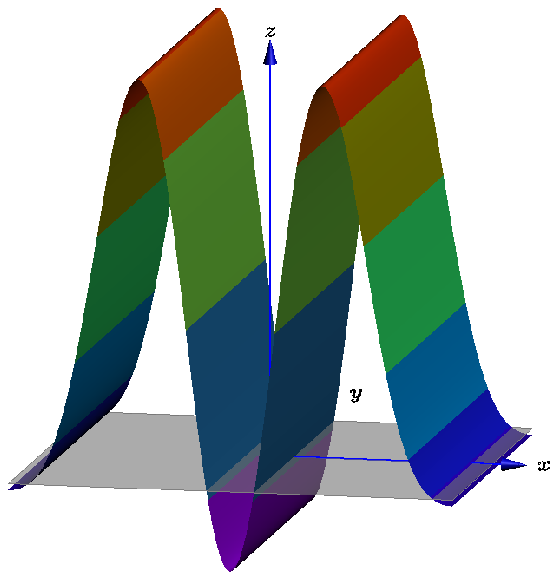}}}
\caption{\label{delta3harmstd}  The graph of the function $\delta$ with respect to the plane $xOy$, associated to  the standard--like map corresponding to a fixed three harmonic potential (see the text)
and respectively  to perturbation parameters $\epsilon=0.05, 0.35, 0.48$}
\end{figure}

  In Fig.\ref{deltastd} we illustrate the graph of $\delta$ over the phase space $\mathbb{A}$ identified to $[-0.5, 0.5)\times \R$ (considered as a subset of the plane $xOy$), for different values of the perturbation parameter $\epsilon$ of the standard map,  while in Fig.\ref{delta3harmstd} is illustrated the graph of  the same function associated to  the standard--like maps  $F_\epsilon$,  defined by a three-harmonic potential,
  \beq\label{ThreeHPotent}
  V(x)=\ds\frac{\eta_1}{2\pi}\cos(2\pi x)+\ds\frac{\eta_2}{4\pi}\cos(4\pi x)+\ds\frac{\eta_3}{6\pi}\cos(6\pi x),\eeq
  \noindent where the parameters $\eta_i$, $i=\overline{1,3}$, are respectively $0.18,  -0.42, -0.11$.\par
We note that in variational setting, the function $\delta(x,y(x,x'))$ is the Hessian of the orbital action $W(x,x')=h(x,x')$,  corresponding to  the $1$--segment of trajectory $(x,y)$.\par
Moreover, the restriction of the function $\delta$ to the points of a $(p,q)$--sequence gives the diagonal entries of the Hessian at that sequence.
Actually, if $(x_i,y_i)$, $i=\ovl{0,q-1}$, are points
on a $(p,q)$--type orbit, and
$$D_{(x_i,y_i)}F=\left(\barr{ll} a_i& b_i\\c_i &d_i\earr\right)$$
then the Hessian expressed in $(x,y)$--coordinates is:

\beq\label{hessianxy} {H}_q=\left[\begin{array}{ccccc} \hat{\alpha}_0&\hat{\beta}_0&0&\ldots&\hat{\beta}_{q-1}\\
\hat{\beta}_0&\hat{\alpha}_1&\hat{\beta}_1&\ldots&0\\
0&\ddots&\ddots&\ddots&0\\
0&\ldots&\hat{\beta}_{q-3}&\hat{\alpha}_{q-2}&\hat{\beta}_{q-2}\\
\hat{\beta}_{q-1}&0&\ldots&\hat{\beta}_{q-2}&\hat{\alpha}_{q-1}\end{array}\right],\eeq
where $\hat{\alpha}_i=d_{i-1}b_{i-1}^{-1}+a_ib_{i}^{-1}$,
$\hat{\beta}_i=-1/b_i$, $i=\ovl{0,q-1}$.\par

 \begin{lemma}\label{lemmaeigendelta}  Let  $(x_i,y_i)$ be a $(p,q)$--type orbit, and $H_q$ its Hessian matrix (\ref{hessianxy}). If $\lam_{min}, \lam_{max}$ are the smallest, respectively the largest eigenvalues of $H_q$ then the values of the $1$--cone function, $\delta(x_i, y_i)$, at each point of the orbit obeys the inequalities:
 \beq\label{ineqeigdiff}
 \lam_{min}\leq \delta(x_i, y_i)\leq \lam_{max}\eeq
 \end{lemma}
 Proof:
 The smallest and the largest
 eigenvalues of $H_q$ are as follows:
 \beq\label{eigendeltarel}
\lam_{min}=\min _{v\neq 0} R(v),\quad \lam_{max}=\max_{v\neq 0}R(v),\eeq
\noindent where $R(v)=\dfrac{<v, H_qv>}{<v,v>}$, $v\neq 0$, is the Rayleigh quotient associated to the symmetric matrix $H_q$.
Because $\lam_{min}\leq R(v)\leq \lam_{max}, \all\: v\neq 0$, and $R(e_i)=\delta(x_i,y_i)$, $\all\: i=\ovl{1,q}$, where $e_i$ is a vector of the standard basis in $\R^q$, it follows the  sequence of inequalities.\par\medskip

By the Gershgorin Theorem \cite{carlmeyer} each eigenvalue, $\lam$, of a symmetric matrix $A=(a_{ij})\in\R^{n\times n}$, belongs at least to one of the intervals
$$[a_{ii}-r_i(A), a_{ii}+r_i(A)], \quad r_i(A)=\sum\limits_{j=1,\: j\neq i}^n |a_{ij}|, \quad i=\overline{1,n}$$
It follows that in the case of the matrix $H_q$, associated to a periodic orbit,  the smallest eigenvalue  belongs to the interval:
$$\lam_{min}\in[\delta_k-r_k(A),\delta_k],\quad \delta_k=\min_{i=\overline{0,q-1}}\delta(x_i,y_i)$$

Since the smallest value of the function $\delta$ restricted to a periodic orbit is an upper bound for the smallest eigenvalue of $H_q$, whose sign is an indicator for the Aubry--Mather periodic orbits, the location of the absolute minimum of the function can give valuable information on extremal type of these orbits.\par

Let us study the influence of  location of minima of  the function $\delta$, on the dynamical behaviour of   the standard--like maps.\par

   We restrict our study to  standard--like maps corresponding to $\epsilon\geq 0$. 
Standard--like maps are \pd--twist APMs, that exhibit  folding regions beyond a threshold $\epsilon'$, which  is the smallest parameter $\epsilon$, for which the set  $\mathcal{F}=\{(x,y)\: | \: a(x,y)\leq 0\}$ is non--empty or equivalently the set $\{x\in[0,1)\:|\: V''(x)\geq 1/\epsilon\} \neq \emptyset$.\par

Standard--like maps also exhibit strong folding property, beyond  a parameter $\epsilon^*>\epsilon'$ which is the smallest   parameter for which
the $1$--cone crossing occurs within the folding region: $\delta(x,y)=2-\epsilon V''(x)\leq 0$. Both folding  and strong folding regions, are union of vertical strips, i.e. subsets of the form  $\{(x,y)\,|\, A\leq x \leq B, y\in\R\}$. \par

Depending on the potential $V$, the corresponding standard--like map can exhibit one or more  strong folding regions.
Since $V$ is an even periodic function, $V''$ is also even, and its critical point set contains the points $x=0, x=0.5$. If at least one of these two points is a maximum point for $V''$, then beyond a perturbation parameter,  the corresponding symmetry line (i.e., $\Gamma_0$ or/and $\Gamma'_0$) is included within the strong folding region.  Such a line is a line of global or local minima for the function $\delta(x,y)$, because  $\delta$ does not depend on $y$ (see Figs. \ref{deltastd}, \ref{delta3harmstd}).  \par

\begin{rmk}\label{remarkdelta} If one of the symmetry lines $\Gamma_0, \Gamma'_0$, of  equation $x=x_0$, $x_0\in\{-0.5, 0\}$,
is a line of minima for the function $\delta$,  then the restriction of the function $\delta$
   to that symmetry line, $\delta(x_0, y;\epsilon)$,  is a decreasing
function with respect to $\epsilon$.\end{rmk}

 We call USF  map (uni--component strong folding region map)  a standard--like map   exhibiting only one connected strong folding region, and TSF (two--component strong folding region map), a map exhibiting a strong folding region with two  connected components (one including the symmetry line $\Gamma_0$, and another, $\Gamma'_0$).\par
  Standard--like maps defined by a potential $V$, whose  second derivative has only two critical points   is an USF  map.
Classical standard map (\ref{stdmapdef}) obeys this condition, as well as  some maps having an analytic potential with infinitely many harmonics, such as the  map:
\beq\label{analyticmharm}
\barr{lll} x'&=&x+y'\\
           y'&=&\ds\frac{\epsilon}{2\pi}\ds\frac{\sin 2\pi x}{1-a\cos 2\pi x}\earr, |a|<1\eeq

For $a=-0.3$, for example, it is  an USF map.\par

The map  corresponding to the three harmonic potential (\ref{ThreeHPotent})
is a TSF map.  In Fig.\ref{delta3harmstd} one can see that the function $\delta$ associated to this map
gets negative beyond a parameter, within two vertical strips.\par
Obviously in the case of a potential given by a trigonometric polynomial the coefficients of the polynomial  can be tuned such that the map be a USF map.


 In Section \ref{TwistOrderP} we  show how  the number  of strong folding components    influences  the dynamical behaviour of a standard--like map, respectively the order properties of periodic orbits of large absolute twist number.\par

In the sequel  we  reconsider the definition of the twist number and derive new properties. We measure the twist number  of an orbit with respect to the canonical trivialization of the tangent bundle of the annulus.
The general framework of its study is
based  on the structure of $\csl2$, and properties of translation number of circles maps given  in the
 \ref{apendix1}, and \ref{apendix2}.\par

\section{The twist number revisited}\label{twistnosect}

Let $f$ be a twist area preserving (APM) of the annulus $\A$, and $z_0, z_1, \ldots,
z_{q-1}$ a $q$--periodic orbit. We identify  the
tangent spaces at $z_i$, $T_{z_{i}}\mathbb{A}$, $i=\ovl{0,q-1}$,
with $\R^2$. If $v=v_0$ is a tangent vector at a point $z$ of the
orbit (let us say that $v\in T_{z_0}\mathbb{A}$), then the
corresponding (forward) tangent orbit is:
$$(z_k, v_k),  \quad
\mbox{where}\,\, v_{k+1}=D_{z_{k}}f (v_k), k\in \mathbb{N}$$

The amount of rotation about a periodic orbit is the average
rotation of tangent vectors under the action of the tangent  map
along that orbit. In order to characterize  the amount of rotation
we interpret the derivative matrices, $D_zf\in SL(2, \R)$, as circle
maps (see  \ref{apendix1}), and a particular lift $\widetilde{D_zf}$
allows then to measure how much is rotated a vector $v\in
T_z\mathbb{A}$, by $D_zf$.\par

To be more precise, we consider  $\mathbb{S}^1$   as being either  the
unit circle, $C$, in $\R^2$ or as the quotient group $\R/\Z$.\par  Let
 $\left[\barr{ll}
a&b\\c&d\earr\right]\in\sl2$ be  the matrix $D_zf$, with $b$ positive
or negative depending on whether $f$ has a positive twist or a
negative twist. $D_zf$ maps the vector $(0,1)^T$ to $(b,d)^T$, i.e  the  twist map $f$ tilts the vertical line through  $z\in\A$ to the right ($b>0$) or to the left ($b<0$). Regarding $D_zf$  as a circle map this means
that the point $1/4+\Z\in\mathbb{S}^1$ (the north pole in $\mathbb{S}^1\equiv C$) is displaced clockwise/anticlockwise  less than $1/2$.\par

 This property of $D_zf$,  transferred to a particular lift, amounts to set a condition on its displacement map value at $x=\ds\frac{1}{4}$ ($x$ is a lift to $\R$ of $\ds\frac{1}{4}+\Z\in\mathbb{S}^1$):

\begin{defn} A lift $\widetilde{D_zf}$ of
the derivative matrix, $D_zf$, of a positive $($negative$)$ twist
map, such that its displacement map satisfies, $\Psi(1/4)\in(-1/2,
0)$ $(\Psi(1/4)\in(0, 1/2))$ is called basic lift. \end{defn}

\begin{figure}[h]
\begin{center}
\includegraphics{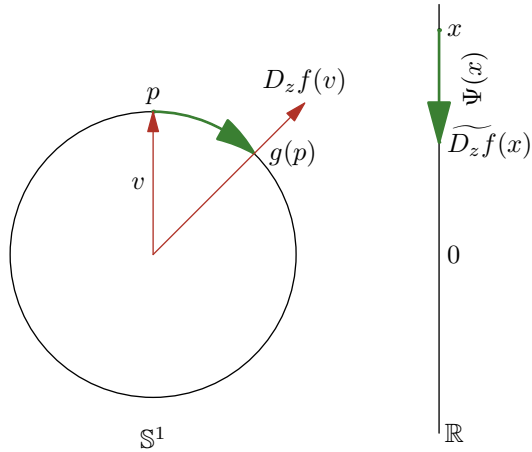}
\caption{\label{Fig0}  Illustration of the action of $D_zf$ on a unit vector $v\equiv p\in \mathbb{S}^1$, and of the corresponding basic lift, on a point $x\in\pi^{-1}(p)$.  $g$ denotes the circle map defined by the matrix $D_zf$.}
\end{center}
\end{figure}

In Fig.\ref{Fig0}  we illustrate the action  of $D_zf$  on  the unit vector $v=(0,1)^T$ or equivalently on $p=  1/4+\Z\in\mathbb{S}^1$, when  $b>0$, as well as how the basic
lift $\widetilde{D_zf}$ displaces the point $x=1/4$, $\pi(x)=p$,  along $\R$ ($\pi:\R\to\mathbb{S}^1$, $\pi(x)= x+\Z$).\par

Let $\widetilde{D_{z_k}f}$ be  the basic lifts associated to
matrices $D_{z_k}f\in SL(2,\R)$, $k=\overline{0,q-1}$, where $z_0,
z_1, \ldots, z_{q-1}\in\mathbb{A}$ are points of a $(p,q)$--periodic
orbit,
 of least period $q$. Because the
covering projection $\Pi:\csl2\to\sl2$ is a group homomorphism, each
lift
$G_k=\widetilde{D_{f^{k-1}(z)}f}\circ\cdots\circ\widetilde{D_{f(z)}f}\circ\widetilde{D_{z}f}$
is a lift of $D_zf^k$, $k=\ovl{1,q}$. We denote by $\Psi_k$ the
corresponding displacement map, and  call $G_k$ the {\it the
$k$--basic lift} associated to the $q$--periodic orbit,
$k=\ovl{1,q}$.\par \bdf The twist number of a $(p,q)$--periodic
orbit, $\mathcal{O}(z)$, is
the translation number of the associated $q$--basic  lift, $G_q$:
\beq\label{defTwistNo}\tau(z)=\lim_{n\to\infty}\displaystyle\frac{(\widetilde{D_{f^{q-1}(z)}f}\circ\cdots\circ\widetilde{D_{f(z)}f}\circ\widetilde{D_{z}f})^n(x)-x}{n}\eeq
\edf

 If
$x\in \R$   then its displacement:
$$\Psi_q(x)=\widetilde{D_{f^{q-1}(z)}f}\circ\cdots\circ\widetilde{D_{f(z)}f}\circ\widetilde{D_{z}f}(x)-x$$
measures the amount of rotation in a period, of a vector, $v\in
T_z\A\equiv\R^2$,  colinear with   $u=(\cos{2\pi x}, \sin{2\pi
x})^T$,  $v=\mu u$, $\mu>0$.\par

The  twist number defined as a translation number, is a
$\csl2$--conjugacy invariant, and thus it does not depend on the
point $z$ of the orbit,  to which one associates the $q$--basic
lift.
\par
In the sequel we discuss the twist number of periodic orbits of
positive twist maps. In this case the basic lift,
$\widetilde{D_zf}$, at an arbitrary point $z\in\mathbb{A}$, belongs
to one of the subsets: $\mathcal{H}_0$, $\mathcal{E}_{-1},
\mathcal{H}'_{-1}\subset \csl2$ or to the common boundary of two of them (see
Fig.\ref{liftsl2r}). Thus the basic lifts $\widetilde{D_{f^k(z)}}$,
$k=\ovl{0,q-1}$, associated to a $q$--periodic orbit have  the
translation number less or equal to zero.\par

The positive twist property of the map leads also to:
\begin{lemma}\label{lemapsi14} If $z$ belongs to a $q$--periodic orbit, and for some
$1<k<q$, the principal lift $G_k$ associated to $z$ has  its displacement at $1/4$,
less than $j/2$, i.e. $\Psi_k(1/4)<j/2$, for some negative integer,
$j$, then no one of the lifts $G_{k+1}, \ldots, G_q$ can have the
displacement at $1/4$, greater than $j/2$.\end{lemma}

In view of this Lemma and taking into account that the twist number
of a $q$--periodic orbit belongs to the range of the displacement
map, $\Psi_q$, i.e. $\tau \in [m(G_q), M(G_q)]$, we can conclude
that the twist number of a periodic orbit of a positive twist map
can be zero or negative.
\par

The twist number defined as above  is minus the twist number defined
by Angenent \cite{angenent}. The negative sign illustrates the
clockwise sense of rotation about the periodic orbit. Its natural
sign allows a better understanding and explanation of its properties
in connection with the displacement map of the lift involved in its
definition.
\par In \cite{mather} the amount of rotation about a periodic orbit
of a positive twist map is defined in a similar way, but in order to
get a positive twist number, it was chosen for each matrix
$D_{f^k(z)}f$, $k=\ovl{0,q-1}$, a lift $G$, with the translation
number $\tau(G)\in[0,1/2)$. Let us argue that such a choice of a
lift is not always possible.\par First we recall shortly the method
used in \cite {mather}  to establish a relationship between the
twist number of a periodic $(p,q)$--type orbit,  and the Morse index
of the corresponding critical sequence of the $W_{pq}$--action.
\par By a non-symplectic  change of coordinates: \beq\label{changecoordxxp}(x,x')\mapsto(x,
y(x,x'))\eeq the derivative $D_zF^q=D_{F^{q-1}(z)}F\cdots
D_{F(z)}F\cdot D_zF\in\sl2$, at a point $z=(x,y)$
of a $(p,q)$--orbit,  is expressed as:

\beq\label{DFzxxp} M=\left[\barr{rr} 0&1\\
-\ds\frac{\beta_{q-1}}{\beta_{0}}&-\ds\frac{\alpha_{0}}{\beta_{0}}\earr\right]\cdots\left[\barr{rr} 0&1\\
-\ds\frac{\beta_{1}}{\beta_{2}}&-\ds\frac{\alpha_{2}}{\beta_{2}}\earr\right]\left[\barr{rr} 0&1\\
-\ds\frac{\beta_{0}}{\beta_{1}}&-\ds\frac{\alpha_{1}}{\beta_{1}}\earr\right]\eeq
The entries $\alpha_i, \beta_i$  are defined by:
\beq\label{alphaibetaih}\beta_i=h_{12}(x_i, x_{i+1}), \quad
\alpha_i=h_{11}(x_i, x_{i+1})+h_{22}(x_{i-1}, x_i),\quad
i\in\Z\,\,(\mbox{mod}\,\, q),\eeq  where $h$ is the Lagrangian
generating function of $F$, and $h_{k\ell}$, $k, \ell\in\{1,2\}$,
are second order partial derivatives of $h$.\par

Thus $D_zF^q$ is $GL^+(2,\R)$--conjugated to the matrix $M$ ($GL^+(2,\R)$ is the group of real $2\times
2$--matrices of positive determinant),  and the
twist number is defined as the translation number of some lift of
$M\in GL^+(2,\R)$ . \par  The lifts $G_B$ of
circle maps, $g_B$, defined by matrices $B\in GL^+(2,\R)$ were
considered in \cite{mather} as having the properties of general
circle maps. The lack of information that the difference between the
maximum and minimum displacement of a lift $G_B$ is less than $1/2$
led to a non concordance.\par

Let us point out that the lifts $G_B$, $B\in GL^+(2,\R)$, have  the
same particularities as the lifts in $\csl2$. Indeed,  to each
matrix $B\in GL^+(2,\R)$ one associates the matrix
$A=\ds\frac{1}{\sqrt{\mbox{det}(B)}}B\in\sl2$, and as a consequence
$GL^+(2,\R)$ can be identified to the direct product of
multiplicative groups, $\sl2\times (0,\infty)$: $$B\in
GL^+(2,\R)\mapsto \left(\ds\frac{1}{\sqrt{\mbox{det}(B)}}B,
\sqrt{\mbox{det}(B)}\,\right)\in \sl2\times (0,\infty)$$ Thus a
 lift of the circle map defined by a matrix in
$B\in GL^+(2,\R)$ is  a pair $(G, \zeta)\in \csl2 \times
(0,\infty)$, such that $B= \zeta\cdot\Pi(G)$, where $\Pi:\csl2\to
\sl2$ is the covering projection. Therefore all properties of lifts
in $\csl2$ and their translation numbers are shared with the lifts
$G_B$, $B\in {GL^+(2,\R)}$.\par Now we are able to show that if
$(x_i, y_i)_{i\in\Z}$ is a $(p,q)$--type orbit of a positive twist
map, $F$, then some matrices $D_{(x_i,y_i)}F\in \sl2$, and hence
their $GL^+(2,\R)$--conjugate matrices that represent the derivative
in $(x,x')$--coordinate, cannot have
 a lift of translation number in $(0,1/2)$, as it was supposed in \cite{mather}. Actually,  if
for some $j$, $D_{(x_j, y_j)}F$ has the trace in  $(-2, 2)$, its
basic lift, $G$, belongs to the subset
 $\mathcal{E}_{-1}\subset \widetilde{SL(2,\R)}$, i.e. $\tau(G)\in(-1/2, 0)$, and any other lift $G+j$, $j\in\Z$, has the
translation number in the interval $(-1/2+j, j)$. Thus no lift $G+j$
can have $\tau\in (0,1/2)$.

Modifying correspondingly the Mather proof in \cite{mather}, that
is taking as lift for each factor in (\ref{DFzxxp})  the basic lift,
one gets the following result:
 \bth
\label{TheorMorseidx} If $(x_n, y_n), n\in\Z$ is a $(p,q)$--type
orbit of the positive twist map $F$, $H_q$ is the Hessian of the
action $W_{pq}$ at ${\bf x}=(x_n)$, $I$ the Morse index of $(x_n)$,
and $\tau$ is the twist number of the corresponding $f$--orbit, then
$I$ and $\tau$ are related as follows:\par
  i) If $I$ is even, then $\tau=-I/2$;\par
  ii) If $I$ is odd, then $\lceil -I/2 \rceil -1\leq \tau < \lceil -I/2 \rceil$, with  equality if and only if {\bf x}
is a degenerate critical point. $\lceil a\rceil$ denotes the least
integer number greater or equal to  $a\in\R$.\eth

From Theorem \ref{TheorMorseidx}, and the general properties of
lifts in $\csl2$ we have the following:

\begin{cor}\label{CorrolaryMI}

i)  The twist number of a minimizing orbit is zero.\par
ii)   Along an elliptic  mini-maximizing  $q$--periodic orbit, the
maps $f^k$, $k=\ovl{1,q}$ have positive twist property.
\end{cor}

Proof:  ii) Actually, if $f^k$, $k=\ovl{2,q}$ violated the twist
property, then the $k$-- basic lift, $G_k$, associated to the orbit
would have a displacement, such that $\Psi_k(1/4)\leq -1/2$,  thus
contradicting, in view of Lemma \ref{lemapsi14},   that the twist
number is $\tau\in (-1/2, 0)$.
\par\medskip

Naturally,  one raises the following questions: under what
conditions a periodic orbit of twist number less than $-1/2$ can
exist? Do such orbits so rarely appear  that no reference to their
existence in the dynamics of most studied twist maps can be found in
literature?\par The Proposition below points out  that such orbits
typically appear at a period doubling bifurcation (details on this type of bifurcation can be found in \cite{mackay93}, \cite{meyerhall}):
\begin{prop}\label{propperiodD} Let $f_\epsilon$, $\epsilon\geq 0$, be a positive twist map, which is a perturbation of an integrable map, $f_0$.
 If at some
threshold $\epsilon'$,  a period doubling bifurcation of a
mini--maximizing $(p,q)$--periodic orbit occurs, with transition
from  elliptic to inverse hyperbolic orbit,  then the new born
$2q$--periodic orbit is elliptic, having the twist number within the
interval $(-1,-1/2)$. If an inverse period doubling bifurcation
occurs, i.e. with transition from an inverse hyperbolic to an
elliptic $q$--periodic orbit, then at the threshold of bifurcation a
new regular hyperbolic $2q$--periodic orbit emerges, and its twist
number is -1.\end{prop}

\begin{figure}[h]
\centerline{\includegraphics{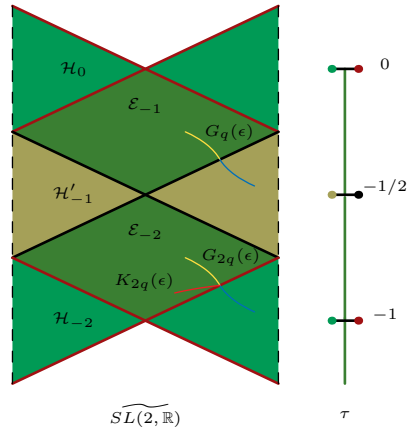}}
 \caption{\label{liftulSL2Rperd} Illustration of the paths in the universal covering space $\csl2$, related to the period doubling bifurcation of a $q$--periodic orbit. }
\end{figure}

Proof: We prove only the first part, the proof of the second one
being similar. Let $\epsilon\mapsto z_{\epsilon}$, $\epsilon \in
(\epsilon_1, \epsilon_2)$, be the continuous path in the annulus
described by the $(p,q)$--periodic orbit that at
$\epsilon'\in(\epsilon_1, \epsilon_2)$ undergoes a period doubling
bifurcation. One associates to this path the
continuous path, $\epsilon\mapsto G_q(\epsilon)$, in the covering
space, $\csl2$, described by the corresponding $q$--basic lift. The
latter one starts in the region $\mathcal{E}_{-1}$ (Fig.
\ref{liftulSL2Rperd}), intersects transversally the parabolic cone
$\mathcal{P}_{2k-1}^{-1+}$, and enters $\mathcal{H}'_{-1}$.\par

On the other hand, each $z_\epsilon$ is also a fixed point of
$f^{2q}_\epsilon$,  and thus the associated $2q$--basic lift is
$G_q^2(\epsilon)$. The twist number along the path $\epsilon \mapsto
G_q^2(\epsilon)$ is double the twist number  along the path
$\epsilon\mapsto G_q(\epsilon)$, because
$\tau(G_q^2(\epsilon))=2\tau(G_q(\epsilon)))$, i.e.
$\tau(G_q^2(\epsilon))\in(-1, 0)$, for $\epsilon \in (\epsilon_1,
\epsilon')$,  it gets $-1$ at $\epsilon'$, and remains $-1$ within
the interval $(\epsilon', \epsilon_2)$. Thus $G_q^2(\epsilon)$
describes  a path from $\mathcal{E}_{-2}$ to $\mathcal{H}_{-2}$,
that intersects transversally the parabolic cone,
$\mathcal{P}_{2k}^{1-}$. From the point of intersection to
$\mathcal{P}_{2k}^{1-}$, a new path emerges, namely the path
described by the $2q$--basic lift $K_{2q}(\epsilon)$, associated to
the new born $2q$--periodic orbit, $p_\epsilon$,
$\epsilon\in(\epsilon', \epsilon_2)$. Because of  the continuity of
the maps $\epsilon\in [\epsilon', \epsilon_2)\to p_\epsilon\to
K_{2q}(\epsilon)$,  and $\tau:\csl2\to\R$, the new path enters the
region $\mathcal{E}_{2}$. Thus the new born elliptic $2q$--periodic
orbit has the twist number within the interval $(-1,
-1/2)$.\par\medskip

Because the $2q$--basic lift associated to the new born elliptic
$2q$--periodic orbit  has the range of its displacement map,
$[m(K_{2q}), M(K_{2q})]\subset (-1, -1/2)$, each tangent vector at a
point of that orbit  is rotated clockwise  in one period,  with an
angle between $\pi$ and $2\pi$.\par

By Theorem \ref{TheorMorseidx} the twist number of a periodic orbit
is well defined by an even Morse index of the corresponding critical
sequence,
 but there is an ambiguity in
deciding the demi-unit interval (an interval of length 1/2) that contains the twist number, when
the Morse index is odd. For instance if the Morse index is 1, the twist number can belong either to $[-1/2, 0)$ or to $(-1, -1/2)$.
\par Because the translation number is  not a group homomorphism
(see \ref{apendix1}), we cannot deduce  the twist number of a
$(p,q)$--periodic orbit, from the translation numbers of the basic
lifts, $\widetilde{D_{z}f}, \ldots, \widetilde{D_{f^{q-1}(z)}f}$,
associated to that orbit.\par

If the Morse index of a $(p,q)$--critical sequence  is an odd number, and the corresponding orbit is elliptic,  we can use a naive method to get the demi--unit interval
(i.e. an interval of the form $(-(k+1)/2, -k/2)$, $k\in\N$) containing its twist number.
 Such a method is based
on the property that if $z, z'$ are any two  points on the same
(elliptic in our case) $(p,q)$--periodic orbit, then $D_zf^q$ and
$D_{z'}f^q$, are $\sl2$--conjugated. Hence   the corresponding
$q$--basic lifts, $G_q$, $H_{q}$, are $\csl2$--conjugated and by
Proposition \ref{conjugminmax}, their  displacement maps  have the range  included in the same demi--unit
interval $(-(k+1)/2, -k/2)$, for some   integer $k\geq 0$. \par
In order to get the actual value of $k$,  we can choose an arbitrary  point $z$ on the
elliptic $q$--periodic orbit,  and compute a $(q+1)$--length segment of
tangent trajectory, $(z_i, w_i)$, $w_{i+1}= D_{z_i}f (w_i)$,
$i=\ovl{0, q}$, with $z_0=z$ and  $w_0=(1,0)^T\in T_{z}\mathbb{A}$.
Let  $c_i$ be  the $x$--coordinates of the vectors $w_i$,
$i=\ovl{0,q}$. If   $k$ changes of sign in the sequence $(c_i)$,
$i=\ovl{0,q}$, are recorded, then  the twist
number belongs to the demi--unit interval, $(-(k+1)/2,
-k/2)$.\par

This method works well  for low values of the period q. In this way
 the reader can check easily the result in Proposition
\ref{propperiodD} in the case of  standard map (Eq.\ref{stdmapdef}) corresponding to
$\varepsilon=2$, respectively $\varepsilon=4$, for which period doubling bifurcation of
the $2$--periodic orbit $(0, 0.5)$, respectively of the fixed point
$(0,0)$, occurs.\par

 The naive method is not suitable for
hyperbolic periodic orbits, because in this case $m(G_q)$ and $M(G_q)$ belong to
distinct demi-unit intervals.\par In the sequel,
 we present  a general  method to deduce the smallest interval $\mathcal{I}$ of ends in $\Z/2$ that contains the interval $[m(G_q), M(G_q)]$
of  the twist number of a $(p,q)$--periodic orbit. The method is
based on the spectral properties of the Hessian matrix, $H_q$, and
of an associated symmetric matrix.
\par

The Hessian  matrix   $H_q$ of the $W_{pq}$--action associated to a
$(p,q)$--periodic orbit, $q\geq 3$, is a Jacobi periodic matrix (i.e. a
symmetric tridiagonal matrix with non-null entries in the upper
right, and left lower corners, and the next--to--diagonal entries
have the same sign). For $q=2$, $H_2$ is simply a symmetric matrix. \beq\label{hessianWpq}
H_q=\left[\begin{array}{ccccc} \alpha_0&\beta_0&0&\ldots&\beta_{q-1}\\
\beta_0&\alpha_1&\beta_1&\ldots&0\\
0&\ddots&\ddots&\ddots&0\\
0&\ldots&\beta_{q-3}&\alpha_{q-2}&\beta_{q-2}\\
\beta_{q-1}&0&\ldots&\beta_{q-2}&\alpha_{q-1}\end{array}\right],
\quad q\geq 3,\quad H_2=\left[\begin{array}{cc} \alpha_0&\beta_0+\beta_1\\ \beta_0+\beta_1& \alpha_1\end{array}\right],\eeq where $\alpha_i\in\R, \beta_i<0$,
$i=\overline{0,q-1}$, are defined in (\ref{alphaibetaih}). Such a
Jacobi matrix is called  periodic because the sequences
$(\alpha_i)$, $(\beta_i)$, $i\in\Z$,  are  $q$--periodic.
\par The spectral properties of periodic Jacobi matrices were
studied extensively in the late  seventies. In his approach to the
twist number, Angenent \cite{angenent}  exploited a result  by van
Moerbeke \cite{moerbeke}. The partial characterization of the
spectrum of periodic Jacobi matrices given by van Moerbeke is
complemented  in \cite{ferguson}, \cite{golub}, and a thoroughly
presentation of  its properties, in connection to discrete Hill
equation, can be found in \cite{toda}.
\par Exploiting the properties of periodic Jacobi matrices  we  are able to visualize the results in Theorem
\ref{TheorMorseidx}, and moreover to enhance it. We give a method to
deduce the smallest interval $\mathcal{I}$ of ends in $\Z/2$, containing the twist
number of a $(p,q)$--periodic orbit, from the position of the real
number zero, with respect to an interval whose ends are eigenvalues
of the Hessian $H_q$, and/or eigenvalues
 of the
 symmetric matrix, $H_q^-$, which is a slight modification of $H_q$:\beq\label{hessianWpqmin}
H_q^-=\left[\begin{array}{ccccc} \alpha_0&\beta_0&0&\ldots&-\beta_{q-1}\\
\beta_0&\alpha_1&\beta_1&\ldots&0\\
0&\ddots&\ddots&\ddots&0\\
0&\ldots&\beta_{q-3}&\alpha_{q-2}&\beta_{q-2}\\
-\beta_{q-1}&0&\ldots&\beta_{q-2}&\alpha_{q-1}\end{array}\right], q\geq 3, \quad H_2^-=\left[\begin{array}{cc} \alpha_0&\beta_0-\beta_1\\ \beta_0-\beta_1& \alpha_1\end{array}\right]\eeq

The Theorem \ref{TheorMorseidx} can be enhanced to:
\bth\label{TheorMorseidxN} If $(x_n, y_n), n\in\Z$ is a
$(p,q)$--type orbit of the positive twist map $F$, $H_q$ is the
Hessian of the action $W_{pq}$ at ${\bf x}=(x_n)$, $I$, the Morse
index of $(x_n)$, $\tau$, the twist number of the corresponding
$f$--orbit, and $I'$ is the number of negative eigenvalues of the
symmetric matrix $H_q^-$, associated to $H_q$, then $I, I'$ and
$\tau$ are related as follows:
\begin{itemize}
  \item[i)] If $I$ is even, then $\tau=-I/2$;\par
  \item[ii)] If $I$ is odd and $H_q$ has at least one zero eigenvalue, then
  $\tau=-(I+1)/2$;
  \item[iii)] If $I$ is odd and $H_q$ has no zero eigenvalue, then:
  \begin{itemize}
  \item[a)] If $I'=I-1$, then $\tau\in(-I/2, -(I-1)/2)$;
  \item[b)] If $I'=I$, then $\tau=-I/2$;
  \item[c)] If $I'=I+1$, then $\tau\in (-(I+1)/2, -I/2)$;
  \end{itemize}

  \end{itemize}
\eth

Proof:

\par In order to relate the spectrum of the Hessian $H_q$ to the
trace of the derivative $D_zF^q$, at a point $z=(x,y)$ of a
$(p,q)$--type orbit,
 one associates to the Jacobi periodic matrix, $J=H_q$,
 the discrete  Hill equation \cite{toda}, $J\xi=\lam \xi, \lam \in\R$, where
 $\xi=(\xi_i)$, is a $q$--periodic sequence in $\R$.
 The coordinate-wise version of this equation is:
 \beq\label{HillDiscr}\beta_{i-1}\xi_{i-1}+(\alpha_{i}-\lam)\xi_i+\beta_i\xi_{i+1}=0,\quad
\lam\in\R\eeq

 Among the  solutions $\xi(\lam)=(\xi_i(\lam))$ of this second-order
difference equation, of particular interest is the Floquet solution,
i.e. the solution ,   that satisfies the condition
$\xi_{i+q}(\lam)=\eta(\lam) \xi_i(\lam)$, with $\xi_0(\lam),
\xi_1(\lam)$ given, and   $\eta(\lam)$  an eigenvalue of the
monodromy matrix: \beq\label{monodromyM} M(\lam)=\left[\barr{cc}
0&1\\-\ds\frac{\beta_{q-1}}{\beta_q}&\ds\frac{\lam-\af_q}{\beta_q}\earr\right]\cdots
\left[\barr{cc}
0&1\\-\ds\frac{\beta_{1}}{\beta_2}&\ds\frac{\lam-\af_2}{\beta_2}\earr\right]\left[\barr{cc}
0&1\\-\ds\frac{\beta_{0}}{\beta_1}&\ds\frac{\lam-\af_1}{\beta_1}\earr\right]\eeq

 The matrix $M(\lam)$ belongs to $\sl2$, while its factors to
$GL^+(2,\R)$.
 We note that for $\lam=0$, $M(0)$ is the derivative $D_zF^q$,
expressed in $(x,x')$--coordinates (see (\ref{DFzxxp})).\par

 The trace of the matrix $M(\lam)$ is a polynomial  of order
$q$, with real coefficients,
\beq\label{HilldiscrF}\mbox{tr}(M(\lam))=(\beta_0\beta_1\cdots\beta_{q-1})^{-1}\lam^q+\cdots\eeq
 called the Hill discriminant \cite{toda}. \par
 The characteristic polynomial of
 a
periodic Jacobi matrix, $J$, defining the discrete Hill equation  is
related to the Hill discriminant by the so called {\it Hill formula}
\cite{ferguson}, \cite{toda}: \beq\label{HillFormula}
\mbox{det}(J-\lam I)=(-1)^q\beta_0\beta_1\ldots
\beta_{q-1}(\mbox{tr}(M(\lam))-2)\eeq

  The history of Hill formula
 and  generalization to multidimensional discrete and
continuous Lagrangian systems can be found  in \cite{bolotin}.
\par

From (\ref{HillFormula}) it follows that  the eigenvalues, $\lam_0,
\lam_1,\ldots, \lam_{q-1}$, of the Hessian, $H_q$, are the roots of
the equation $\mbox{tr}(M(\lam))=2$. The eigenvalues  can be ordered
in the following way \cite{toda}: \beq\label{eigenHq}\lam_{2i}<
\lam_{2i+1}\leq\lam_{2i+2}, \all\,\, i\geq 0\eeq

On the other hand the roots of the equation $\mbox{tr}(M(\lam))=-2$
can be ordered \cite{toda} as: \beq\label{eigenHqm} \lam'_{2i}\leq
\lam'_{2i+1}<\lam'_{2i+2}, \all\,\, i\geq 0,\eeq and the two series
are interlaced: \beq\label{allroots} \lam_0<\lam'_0\leq
\lam'_1<\lam_1\leq\lam_2<\lam'_2\leq\cdots \eeq  Being roots of the
equations  $\mbox{tr}(M(\lam))=2$, $\mbox{tr}(M(\lam))=-2$, with
$M(\lam)\in\sl2$, $\lam_{2i+1}$, respectively $\lam'_{2i}$, $i\geq
1$, can have at most double multiplicity. Thus each  root occurs in
(\ref{allroots}) as many times as its multiplicity.\par

For practical purposes it is important to  note that the roots
$\lam'_i$, $i=\ovl{0,q-1}$, are the eigenvalues of the symmetric
matrix $H_q^-$, defined in (\ref{hessianWpqmin}), associated to the
Hessian $H_q$. Indeed, by \cite{golub} the characteristic polynomial
of $H_q^-$ is for any $q\geq 2$:
$$
\mbox{det}(H_q^- -\lam I_q)=\mbox{det}(H_q-\lam I_q) + 4
(-1)^q\beta_0\beta_1\cdots \beta_{q-1}$$
 Comparing to the expression of the Hill discriminant
 (\ref{HilldiscrF}) one gets, that $\mbox{tr}(M(\lam))+2=0$
 iff $\mbox{det}(H_q^- -\lam I_q)=0$.\par
\begin{figure}
\begin{center}
\includegraphics[scale=1.25]{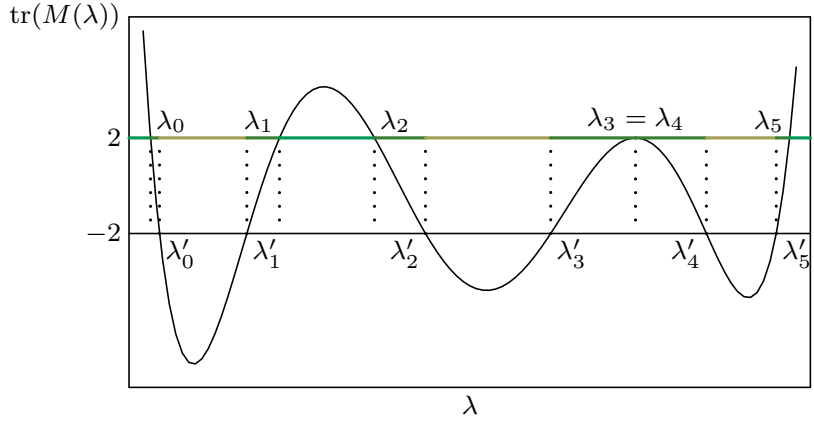}
\caption{\label{Hilldiscr} The graph of the trace of monodromy
matrix, $M(\lam)$,  the eigenvalues of a Hessian matrix $H_6$ and
its associated matrix $H_6^-$.}
\end{center}
\end{figure}

For the reader convenience we illustrate in Fig.\ref{Hilldiscr}  the
graph of a Hill discriminant corresponding to an even integer,
$q$.\par

  Because the leading
coefficient of the polynomial $\mbox{tr}(M(\lam))$ is
$(\beta_0\beta_1\ldots\beta_{q-1})^{-1}$, it follows that in both
cases, $q$ even and odd,
 $\lim_{\lam\to-\infty}\mbox{tr}(M(\lam))=+\infty$. Hence  the  dynamical type of the monodromy matrix
 (Fig.\ref{Hilldiscr}) is as follows:\par
a) If   $\lam<\lam_0$ or $\lam\in (\lam_{2i+1}, \lam_{2i+2})$,
$i\geq 0$ with $\lam_{2i+1}\neq \lam_{2i+2}$, then  $M(\lam)$ is
regular hyperbolic.\par b) For $\lam\in(\lam_{2i}, \lam'_{2i})$ or
$\lam\in(\lam'_{2i+1}, \lam_{2i+1})$, $i\geq 0$, $M(\lam)$  is
elliptic.\par c) For $\lam\in(\lam'_{2i}, \lam'_{2i+1})$,
$\lam'_{2i}\neq\lam'_{2i+1}$, $i\geq 0$, $M(\lam)$ is inverse
hyperbolic.\par\medskip

Now the location of the real number zero with respect to the
interlaced series of eigenvalues (\ref{allroots}) gives information
on the Morse index of $H_q$, the signature of the matrix $H_q^-$, and
 the dynamical type (regular hyperbolic, elliptic or inverse
hyperbolic) of the $(p,q)$--periodic orbit we are studying.\par If
the Morse index $I$ of the Hessian $H_q$ is odd, and $H_q$ has no
zero eigenvalue, i.e. $0\in (\lam_{2i}, \lam_{2i+1})$, for some
$i\geq 0$, then by Theorem \ref{TheorMorseidx}  the twist number
 belongs to an interval of length 1, $\lceil -I/2 \rceil -1< \tau <
\lceil -I/2 \rceil$. It follows that the $q$--basic lift $G_q$
associated to our orbit can belong to $\mathcal{E}_{-2I+1},
\mathcal{H}'_{-2I+1}$, or to $\mathcal{E}_{-2I-1}\subset\csl2$, and
the twist number can belong to different subintervals of the above
1--length interval. In order to deduce the right demi--unit interval
to which the twist number belongs we analyze the position of $0$
with respect to the eigenvalues $\lam'_{2i}, \lam'_{2i+1}$, and
deduce the index $I'$ of the symmetric matrix $H_q^-$. Due to the
continuity of the translation number $\tau:\csl2\to\R$, the case
iii) from the Theorem follows directly.\par Similarly one can prove
the case ii).

\begin{cor}
  The smallest  twist number a  $(p,q)$--periodic orbit of a positive twist APM can have is  $-q/2$.\end{cor}

This property was  also stated in \cite{angenent} but no arguments were given.

 In numerical investigations of the twist number, if the Lagrangian generating function of the twist map is defined
only locally (as in the case of tokamap \cite{balescu}), then one can reduce the computations and improve
accuracy, evaluating the Hessian $H_q$ associated to a periodic
orbit, not in $(x,x')$--coordinates, but  in the original symplectic $(x,y)$--
coordinates. In this case the entries of $H_q$ are as in (\ref{hessianxy}).\par

Thus, a simple algorithm that in the first step calculates
 the eigenvalues of the Hessian (\ref{hessianxy}), order
them increasingly, and then tests the cases described in Theorem
\ref {TheorMorseidxN} (computing also the eigenvalues of $H_q^-$, if
necessary, and ordering them increasingly), returns the  value of
the twist number or the demi--unit interval containing the twist
number.

\par

We note that the $2q$--periodic orbits born through a
period--doubling bifurcation are ordered orbits (Birkhoff orbits).
In the next section we give examples  of badly ordered (non--Birkhoff)
$(p,q)$--periodic orbits that are  born through a saddle--center
bifurcation and have the absolute  twist
number greater than  $1/2$. We also show that the so called two--component strong folding region maps can exhibit ordered periodic orbits of large absolute twist  number.
\section{Twist number and order properties of periodic orbits}\label{TwistOrderP}

Computer experiments reveal that plenty of  pairs of periodic orbits of absolute twist number,  $|\tau|> 1/2$,
can be found in the folding region  of the standard map corresponding to the perturbation parameter  $\varepsilon>1$.  \par
The standard map we consider here is defined by the lift  $F:(x,y)\mapsto(x', y')$:
\beq\label{stdmapdef}
\begin{array}{lll}
x'&=& x+y'\\
y'&=& y-\ds\frac{\varepsilon}{2\pi}\sin(2\pi x), \quad \varepsilon\in\R,
\end{array}
\eeq

\begin{figure}
 \centerline{\scalebox{1.5}{\includegraphics{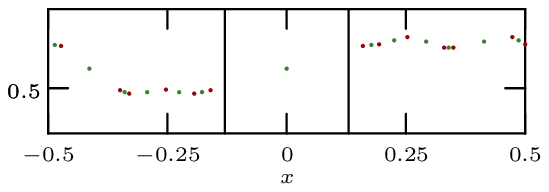}\includegraphics{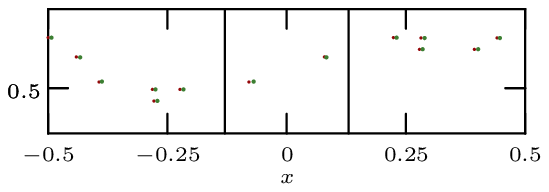}}}
\caption{\label{std813} Four $(8,13)$--periodic orbits of the standard
map corresponding to $\varepsilon=1.4578$. At left is illustrated the pair of
orbits given by Aubry--Mather theory, and at right  a pair of orbits
of twist
number $-1$, respectively  $\tau\in (-1, -1/2)$.}
\end{figure}

{\bf Example 1.} The standard map corresponding to the perturbation
parameter $\varepsilon = 1.4578$ has the minimizing $(8,3)$--periodic orbit represented by $( 0.5,
 0.693821664066481)$, and the  mini--maximizing one by  $(0,
 0.586319735666097)$ (Fig. \ref{std813}, left panel).\par

The point $z=(0.5, 0.7237177352891645)$ is also a $(8,13)$--periodic
point. Applying   the ubiquitous method used in the
study of twist area preserving maps, namely, computing only the
residue, $R=0.44794945246558$, we can mistakenly conclude that it is
a mini--maximizing periodic orbit. Computing  further the
eigenvalues of the Hessian $H_{13}$ we get  the Morse index
$I=1$, and analogously we are tempted to conclude that it is
mini-maximizing (the index equal to one is only a necessary
condition to be mini-maximizing, not a sufficient one). But the orbit is unordered and thus it cannot be mini-maximizing.
 The eigenvalues of the associated matrix
${H}_{13}^-$ reveal  that its twist number belongs to $(-1, -1/2)$.
Hence this orbit
 is an elliptic periodic orbit  with the property that each vector
$v\in T_z\A$ is rotated clockwise about the orbit of $z$, more than
$\pi$ radians in one period.\par

The orbit of $z$ was born through a saddle--center bifurcation at a
slightly lower value of $\varepsilon$, and its pair is the  regular hyperbolic
orbit of the point $z'=(0.5, 0.723372427461496)$. $z'$ has the
twist number  $\tau=-1$, and its  orbit is also  unordered. They are
illustrated in Fig. \ref{std813}, right panel.\par
 \par

One can observe that projecting each of the four orbits onto the first factor, $\mathbb{S}^1$, identified with  $[-0.5, 0.5)$, and  denoting by $x_\ell$, $x_r$ the abscissas of the minimizing--orbit, nearest  to $0$ ($x_\ell$ at left and $x_r$ at right of $0$), then within
the interval $(x_\ell, x_r)$ each projected non--Birkhoff orbit has two points.\par
{\bf Example 2.} The  point
$z=(-0.387429398213243,  0.225141203573513)$ of the  standard map corresponding to $\varepsilon=3$ has  an inverse hyperbolic $(1,7)$--periodic orbit of translation number $\tau=-1.5$,
and  its pair $z'=(0.389952663038891,  0.536877082952284)$  has a  regular hyperbolic $(1,7)$--periodic orbit of twist number $\tau=-1$. \par
For the same parameter the Aubry--Mather orbits of type $(1,7)$ are the following:
$(0,  0.390375216334041475)$, mini--maximizing, and
 $( -0.5,  0.008258333079127612)$, minimizing.\par
 Considering again the interval $(x_\ell, x_r)$ as in the previous example, we note that within this interval each orbit of large absolute twist number
 has three points.\par

\begin{figure}

 \centerline{\scalebox{1.3}{\includegraphics{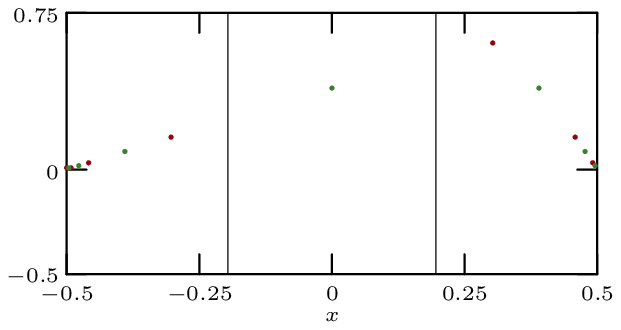}\includegraphics{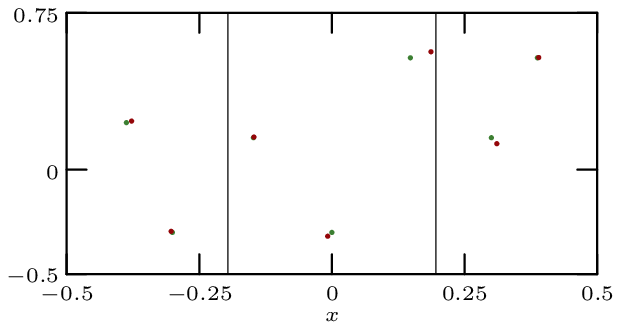}}}
\caption{\label{dynstdmapTW} Four $(1,7)$--periodic orbits of the standard
map corresponding to $\varepsilon=3$. At left is illustrated the pair of
orbits given by Aubry--Mather theory, and at right  a pair of orbits
having the twist number $-1$, respectively $-1.5$.
 }
\end{figure}

In Fig.\ref{dynstdmapTW} are illustrated all four orbits of type $(1,7)$. 

\par
In an attempt to asses numerically the efficiency of the Boyland converse KAM criterion \cite{boyland2}, Leage and MacKay \cite{leage}
looked for badly ordered periodic orbits of the standard map. They detected such orbits for  large perturbations parameters.\par

  We note that the badly ordered  orbits reported in \cite{leage} have also non--zero twist number and  $|\tau|> 1/2$. They are orbits born through a saddle-center bifurcation, at a parameter value greater than the threshold of the   period doubling bifurcation of the mini--maximizing periodic orbit of the same type $(p,q)$.\par
  {\bf Example 3.}  The point $z=(0.181826060531142,  0.681826060531142)$ found in \cite{leage} for  $\varepsilon=7.221365$ is an elliptic  point of type $(1,3)$, twist number
   $\tau\in (-1, -1/2)$, and its pair $z'=(0.181803685309653,  0.681803685309654)$ has a regular hyperbolic $(1,3)$--orbit,  and   $\tau=-1$.\par
   The above examples can suggest that $(p,q)$--periodic orbits of large absolute twist number are  unordered.
A natural question is whether a positive twist map can also exhibit ordered  $(p,q)$--periodic orbit ($p,q$ relative prime integers) of twist number less than $-1/2$. The following example shows that the answer is positive.\par

\begin{figure}
\centerline{\includegraphics{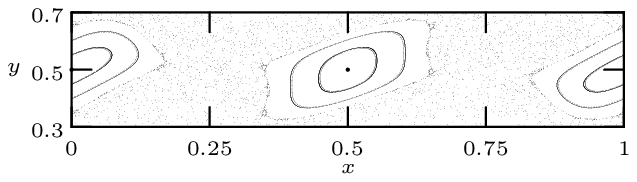}\includegraphics{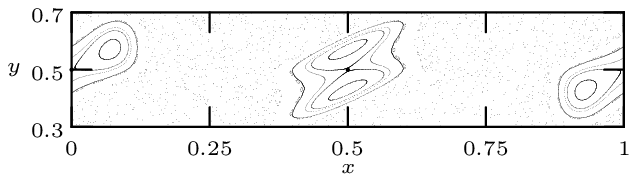}}
\centerline{\includegraphics{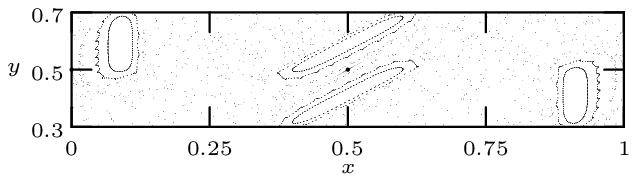}\includegraphics{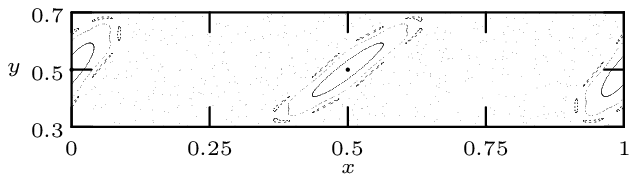}}
\centerline{\includegraphics{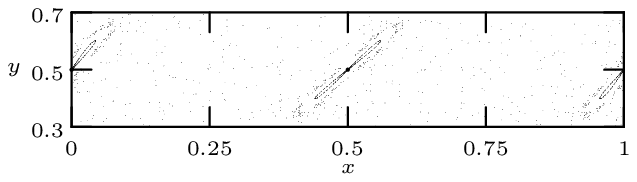}\includegraphics{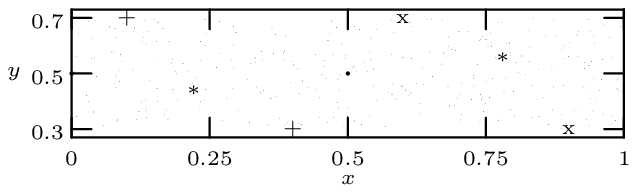}}
\caption{\label{seqbifper} The phase portrait of the map $F_\epsilon$, Eq. (\ref{mapex5eps}),  illustrating the sequence of bifurcations that the orbit of the point $(0,0.5)$ undergoes as $\epsilon$ increases from $0$ to $1$. First line: at left, $\epsilon=0.028$, $(0,0.5)$ has an elliptic mini--maximizing  orbit; at right,  $\epsilon=0.4$, slightly above a  period doubling bifurcation. Second line: at left, $\epsilon=0.46$, slightly above the second period doubling bifurcation (from $\tau=-1/2$ to $\tau\in(-1,-1/2)$); at right, $\epsilon=0.65$, $(0,0.5)$ has an elliptic orbit of twist number $\tau\in (-1, -1/2)$. Third line: at left $\epsilon=0.8$, after a  Rimmer bifurcation (from $\tau\in(-1,-1/2)$ to $\tau=-1$); at right, $\epsilon=1$, $(0,0.5)$ has regular hyperbolic orbit of twist number $\tau=-1$, and the periodic orbits born through Rimmer bifurcation (the orbits marked with $+$, respectively $\rm x$, are inverse hyperbolic of $\tau=-1/2$).}
\end{figure}

{\bf Example 4.}
We consider a multi--harmonic standard map defined by the
 following lift:
 \beq\label{mapex5}
 \barr{lll} x'&=&x+y+ P(x)\\
 y'&=&y+ P(x),\earr\eeq
 \noindent where $P$ is the trigonometric polynomial, $P(x)=\eta_1 \sin(2 \pi  x)+\eta_2\sin (4\pi x)+\eta_3\sin(6\pi x)$, corresponding to  the point
  $\eta=(0.18,  -0.42, -0.11)$ in the parameter space.  This map is reversible  and
 has four $(1,2)$--periodic orbits (Fig.\ref{seqbifper}, third line, right) :\par
a) $(0.220481551875563,  0.440963103751125)$ is a symmetric  minimizing orbit.\par
b) $(0,0.5)$ is a  regular hyperbolic point having the twist number equal to $-1$ (the associated Hessian has  negative eigenvalues $\lambda_1=-5.48881762834532, \lambda_2=-1.06693368771638$,
          and thus the corresponding $(1,2)$--sequence is a sequence of maximum action). \par
c) $z_1=(0.099868507451075,0.699123941758688)$, \newline\indent \hspace{0.28cm} $z_2=(0.900131492548925, 0.300876058241312)$\newline are two inverse hyperbolic points having distinct
 orbits, of twist number $-1/2$. \par
All four $(1,2)$-- periodic orbits are cyclically ordered.\par

We are wondering why in this case an orbit of large absolute twist is ordered, while in the case of standard map such orbits are badly  ordered.
In order to give an answer we connect this map to an integrable one,
i.e. we consider the maps $F_\epsilon$, $\epsilon\in[0,1]$ defined by:
\beq\label{mapex5eps}
\barr{lll} x'&=&x+y+ \epsilon P(x)\\
 y'&=&y+ \epsilon P(x),\earr\eeq

\noindent and we analyze the dynamical type and twist number of $(1,2)$--periodic orbits, as $\epsilon$ increases from $0$ to $1$. \par

\begin{figure}
 \begin{center}
\includegraphics[scale=0.7]{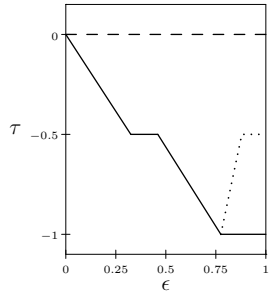}
\caption{\label{twistevol} The twist number of the   $(1,2)$--periodic orbits map (Eq. \ref{mapex5eps}),  as a function of $\epsilon$. The dashed line is the twist number of the minimizing orbit, the continuous line is the twist number of the point $(0, 0.5)$, and the dotted one is the twist number for the two asymmetric periodic orbits, born through a  Rimmer bifurcation.}
\end{center}
\end{figure}

 The map   $F_\epsilon$,  close to the integrable twist map, has   two $(1,2)$--periodic orbits:  a minimizing one that persists down to $\epsilon=1$, and the orbit of the point $(0,0.5)$ which starts as a mini--maximizing orbit of elliptic type (Fig.\ref{seqbifper}, first line, left). The latter one passes through a sequence of bifurcations (Fig.\ref{seqbifper}) that  lead to the following change of its twist number (see Fig.\ref{twistevol}):\par
 elliptic mini--maximizing $\to$ inverse hyperbolic of $\tau=-1/2$ $\to$ elliptic of twist number in $(-1, -1/2)$ $\to$  regular hyperbolic of $\tau=-1$. \par  The first two bifurcations are period doubling, while the last one is a Rimmer--type bifurcation that leads to the creation of the two distinct asymmetric periodic orbits, $\mathcal{O}(z_1)$, $\mathcal{O}(z_2)$.\par
Hence the maximizing orbit $(0, 0.5)$ of the map $F_1$ was born as an ordered mini--maximizing orbit, and through a sequence of bifurcations it changed  the twist number as it is illustrated in Fig.\ref{twistevol}, but not the order property, because  the period doubling, and Rimmer bifurcation of ordered orbits leads to ordered orbits.
\par
We also note that while the $(1,2)$--minimizing orbit does not undergo any bifurcation as the perturbation parameter increases, the (1,2)--mini-maximizing orbit bifurcates, changing also the twist number value after  each period--doubling bifurcation. Moreover at the Rimmer bifurcation threshold, the symmetric mini--maximizing orbit is continued by a maximizing one, and two asymmetric mini--maximizing orbits are born. Thus at $\epsilon=1$ there exist 3 orbits given by Aubry--Mather theory: a symmetric minimizing orbit, and two asymmetric mini--maximizing orbits.\par


\begin{rmk}  From the analysis of the map \ref{mapex5eps} it follows that the mini--maximizing $(p,q)$--periodic orbit of a twist map can have
the twist number $\tau$ either in the interval $[-1/2,0)$ or in the interval $(-1,-1/2)$, while the twist number of the corresponding minimizing orbit is constant and equal to $0$.\end{rmk}\par\medskip
In the sequel we show that the unexpected behaviour of the map (\ref{mapex5eps}) is due to the fact that it is a TSF map. \par
First we show that if a twist map is either USF or TSF we can predict the bifurcations that the orbits intersecting $\Gamma_0^*$ can undergo, as the perturbation parameter varies.
\bpr\label{PropUSFbeh}  Let $F_\epsilon$ be an USF map for $\epsilon\geq\epsilon^*$, and $\Gamma^*$ the symmetry line that get included in   the strong folding region for $\epsilon>\epsilon^*$. If after a slight perturbation of the integrable map $F_0$, a $(p,q)$--periodic orbit, $q>1$,  has a point $(x_0(\epsilon),y_0(\epsilon))$  on $\Gamma^*$, then the orbit can be either mini--maximizing or minimizing.  Increasing $\epsilon$ from  $0$, the orbit can undergo the following bifurcations:\par
a)  If $(x_0(\epsilon),y_0(\epsilon))$ is mini--maximizing of elliptic type  it undergoes necessarily  a period doubling bifurcation at some $0<\epsilon_d\leq \epsilon^*$.\par
b) If $(x_0(\epsilon),y_0(\epsilon))$  is minimizing, then at a threshold $\epsilon_r$ it gets mini--maximizing through a Rimmer--type bifurcation, and at some
$\epsilon_d>\epsilon_r$ it undergoes a period doubling bifurcation, and gets inverse hyperbolic of twist number $\tau=-1/2$.\par

  \epr

Proof.

At  $\epsilon=0$   $\delta(x, y)=2$, for every point  $(x,y)$. Let $H_{q}$ be the Hessian matrix associated to the $(p,q)$--orbit,  having a point on $\Gamma^*$, $H_{q}^-$  the symmetric matrix (\ref{hessianWpqmin}),
and $\lambda_0$, respectively $\lambda'_0$,  the  smallest eigenvalues of these matrices.  From Remark \ref{remarkdelta} it follows that the function $\delta(x_0,y_0; \epsilon)$ decreases from $2$, as $\epsilon$ increases. By Lemma \ref{lemmaeigendelta}  we have that
 $\lambda_0, \lambda'_0\leq \delta(x_0, y_0)<2$, and by Gershgorin Theorem  $\lambda_0$ belongs necessarily to the  interval  $[\delta(x_0, y_0;\epsilon)-2, \delta(x_0, y_0;\epsilon)]$.  Hence  for $\epsilon$ very close to zero, $\lam_0$ can be either negative or positive, i.e.  after a slight perturbation of the integrable map the considered  periodic orbit can be mini--maximizing of elliptic type or minimizing.
 \par  a) If the orbit is mini--maximizing, then since the function $\delta(x_0, y_0;\epsilon)$ is decreasing, the eigenvalue  $\lambda'_0(\epsilon)$  is pushed to the left on the real axis, and at a threshold $\epsilon_d\leq \epsilon^*$ it   crosses $0$,  i.e. the twist number crosses $-1/2$ and a  period doubling bifurcation occurs. \par
 b) If the orbit is minimizing, using similar arguments it follows that necessarily the stated sequence of bifurcations must occur.\par\medskip

\begin{cor} An USF map exhibits a total or a partial dominant symmetry line, that is either for each perturbation parameter or for $\epsilon$ greater than a threshold, $\epsilon^*$, any periodic orbit   having a point on $\Gamma_0^*$ has the twist number less than zero.\end{cor}

It is conjectured that the standard  map has $\Gamma_0$ as a  total dominant symmetry line (as far as we know, no proof was given for this numerical observation).  The map
(Eq.\ref{analyticmharm}) corresponding to $a=-3$ has $\Gamma_0$ as a partial dominant symmetry line.\par


In the case of a TSF map a similar behaviour have the orbits intersecting $\Gamma_0$, but moreover a second period doubling bifurcation from $\tau=-1/2$ to $\tau\in(-1,-1/2)$ is also possible.\par  Let us analyze  the $(1,2)$--periodic orbit of the twist map (\ref{mapex5eps}), intersecting $\Gamma_0$.\par

 After a slight perturbation of the integrable map  the mini--maximizing $(1,2)$--periodic orbit has its  points $(x_0, y_0)$, $(x_1, y_1)$ on the symmetry lines $\Gamma_0, \Gamma'_0$.   The restriction of the $1$--cone function $\delta$ to each of these symmetry lines is a decreasing function with respect to $\epsilon$, namely it decreases from the value $2$ corresponding to $\epsilon=0$, to negative values corresponding to $\epsilon$ for which  $\rm{Cone}_1(\epsilon)$  has two connected components.  It follows  that as $\epsilon$ increases, the smallest eigenvalue, $\lam_0\leq \min\{\delta(x_0, y_0), \delta(x_1, y_1)\}$, of $H_2$ is pushed to the left on the real axis, as well as the smallest eigenvalue $\lam'_0=\min\{\delta(x_0, y_0), \delta(x_1, y_1)\}$ of the matrix $H_2^-$ (see \ref{hessianWpqmin}).  Hence a sequence of bifurcations of the ordered $(1,2)$--periodic orbit occurs leading to a continuous decrease of the twist number. Beyond an $\epsilon$ both points of the orbit are within the strong folding region. In this case
  the Hessian $H_2$ and its associated matrix $H_2^-$ have  negative diagonal entries, and as a consequence
$H_2^-$ has both  eigenvalues negative. Thus the $(1,2)$-periodic orbits intersecting Fix(R) is in this case an ordered orbit, of twist number $\tau<-1/2$. \par

\par
Let us show that  in the case of an USF map  the second period doubling bifurcation  of a mini--maximizing $(1,2)$--periodic orbit is not possible (i.e. transition from  $\tau=-1/2$ to  $\tau\in(-1, -1/2)$).\par

By Proposition \ref{PropUSFbeh} The points of the mini--maximizing $(1,2)$--orbit of an USF map are forced (eventually after a threshold) to belong to the symmetry lines $\Gamma_0, \Gamma'_0$. The points,
$(x'_0, y'_0)$, $(x'_1, y'_1)$, of the minimizing  periodic orbit belong to $\Gamma_1$, respectively $\Gamma'_1$, and  $R(x'_0),y'_0)=(x'_1, y'_1)$. Taking into account that  for standard--like maps  the $1$--cone function is $R$--invariant, i.e. $$\delta(R(x,y))=\delta(x,y),$$ for any point $(x,y)$ on a symmetric periodic orbit, we have the following  Hessian matrices (\ref{hessianWpq})  $H_2$, ${H'}_2$,  associated to the two orbits, and the corresponding symmetric matrices (\ref{hessianWpqmin})  $H_2^-$, ${H'}_2^-$:
\beq\label{Hess2M}
H_2=\left[\barr{cc} \delta(x_0,y_0)&-2\\
-2&\delta(x_1,y_1)\earr\right],\quad       H_2^{-}=\left[\barr{cc} \delta(x_0,y_0)&0\\
0&\delta(x_1,y_1)\earr\right]\eeq
\beq \label{Hess2m}H'_2= \left[\barr{cc} \delta(x'_0,y'_0)&-2\\
-2&\delta(x'_0,y'_0)\earr\right],\quad       {H'}_2^{-}=\left[\barr{cc} \delta(x'_0,y'_0)&0\\
0&\delta(x'_0,y'_0)\earr\right]\eeq

\begin{figure}
\begin{center}
\includegraphics{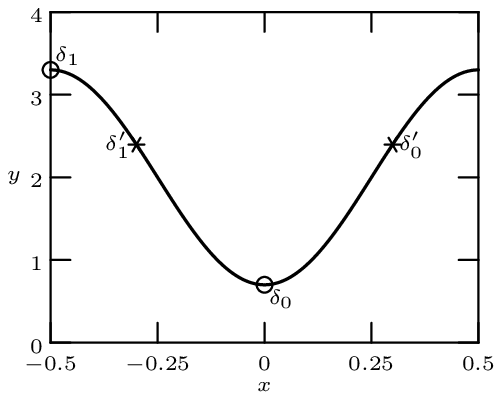}\hspace{1cm}\includegraphics{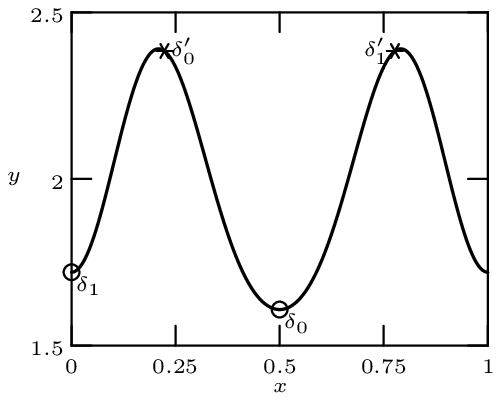}
\caption{\label{delta2Mm} The profile of $\delta$ function associated to the standard map (left), respectively a TSF map (right) and the value of $\delta$ at mini--maximizing points of the $(1,2)$--orbit (marked by $\circ$), and minimizing points (marked by $*$). }
\end{center}
\end{figure}

In Fig.\ref{delta2Mm} is illustrated the profile  of the $\delta$ function for the  standard map, as an USF map, and for the  TSF map (\ref{mapex5eps}), as well as the values of $\delta$ along the mini--maximizing, respectively minimizing $(1,2)$--periodic points of these maps.

\begin{figure}
\begin{center}
\includegraphics{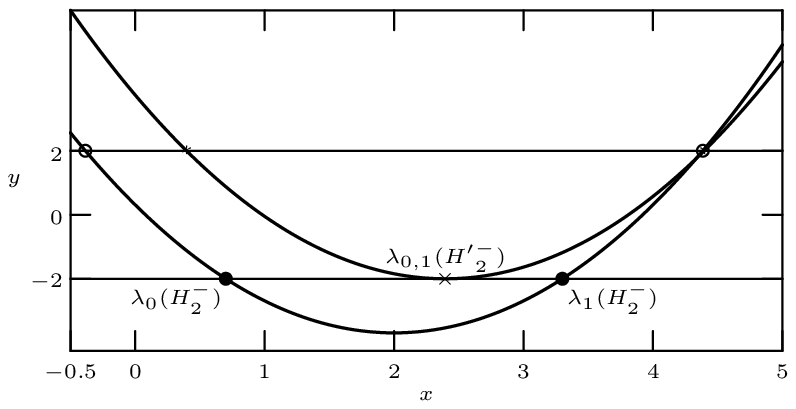}\includegraphics{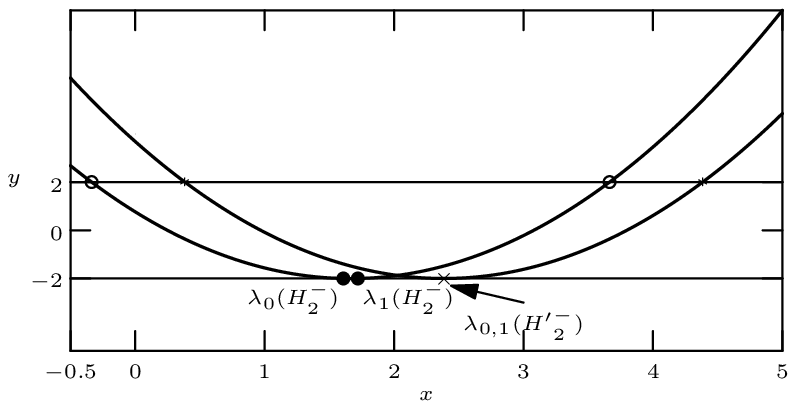}
\caption{\label{poly2TSFUSF} The graphs of the Hill discriminants associated to $(1,2)$--periodic orbits of the standard map (as a USF map), corresponding to $\epsilon=1.3$ (left)
 respectively of the TSF map (\ref{mapex5eps}), corresponding to $\epsilon=0.8$ (right). These maps have the profile of $\delta$ function illustrated in Fig. \ref{delta2Mm}.}
\end{center}
\end{figure}

A simple algebraic computation reveals that the different relative position of diagonal entries of ${H'}_2^-$ with respect to those of $H_2^-$, in the two cases, leads to distinct
relative position of the corresponding eigenvalues  (Fig.\ref{poly2TSFUSF}).  In the case of the  USF map (Fig.\ref{poly2TSFUSF}, left)
the period doubling bifurcation from $\tau=-1/2$ to $\tau\in(-1, -1/2)$  is obstructed because the double eigenvalue of ${H'}_2^-$ is positive, and
less than the second eigenvalue,  $\lambda_1(H_2^-)$, of $H_2^-$.\par
In the case of the TSF map (Fig.\ref{poly2TSFUSF}, right) as the perturbation parameter increases, the double eigenvalues of ${H'}_2^-$ increases, while  $\lambda_0(H_2^-)$, $\lambda_1(H_2^-)$ decrease and  successively they get zero, without any obstruction.\par

We conjecture that in the case of an USF map no   mini--maximizing $(p,q)$--type orbit of twist number $-1/2$ can bifurcate to a an orbit of twist number
$\tau\in(-1,-1/2)$, because of a similar obstruction. In this case  the only periodic orbits of large absolute twist number are born through a saddle--center bifurcation.\par

After the discussion concerning the residue and the index of the  Hessian matrix associated to periodic orbits given in Example 1 we are tempted to  state that
 among residue, the Morse index, and the
twist number of a $(p,q)$--periodic orbit, the twist number  is the
only numerical quantity that leads to a correct decision whether or
not that orbit is one given by Aubry-Mather theory.
However this is not the case, because as the following example illustrates, there can exist badly ordered periodic orbits of twist number zero.\par

{\bf Example 5.} The point $z=( 0.5,  0.323020669964897189)$ of the standard map corresponding to $\varepsilon=2.87$ is a $(3,5)$--regular hyperbolic  point, because its Hessian $H_5$ has the eigenvalues:
$$\barr{l}
0.04121176201302,\quad
   0.53521014626564,\quad
   3.88911846604733,\\
   5.36539881978110,\quad
   5.89051890157580\earr$$
   Its twist number is zero,
   but it is not minimizing because it is badly ordered. \par The Hessian associated to its pair, $z'=( 0.4941763525132936290,  0.3344329777751245200)$, has the eigenvalues:
   $$ \barr{l} -0.07386797172616,\quad
   0.56782907421433,\quad
   3.88372973742859,\\
   5.35939921395639,\quad
   5.88621274711297\earr $$
   The orbit of $z'$  is also badly ordered. The minimizing and mini--maximizing periodic orbits of the type $(3,5)$ are respectively:
   $$(0.5,  0.792447770847267734),\quad  (0,  0.567936349180897305)$$

   The action along the orbits of these four periodic points of type $(3,5)$, in the order of their mention, is respectively:
   $$ 0.891078196747568, \quad 0.891083047652180, \quad   0.804181852816563,\quad 0.831286899496951$$
   \par
   \section{Conclusions}  Deciphering   the structure of the universal covering space of the group $SL(2,\R)$, and exploiting the properties of the translation number of circle maps, we developed a general framework for the  study of the twist number of periodic orbits of twist maps.\par
   We developed an algorithm to compute the value of the twist number or an interval containing this value. Applying this algorithm we illustrated
   that the classical method used in dynamical systems theory, to classify   periodic orbits of a twist map only by their residue or the Morse index of the associated Hessian matrix can lead to errors.\par
We proved that ordered periodic orbits of large absolute twist number are born through a period doubling bifurcation. We also gave examples of periodic orbits born through a saddle--center bifurcation that have large absolute twist number and are unordered.

\par

In order to explain the variation of the twist number in a subclass of standard--like maps we introduced a new tool of study of the dynamical behaviour of these maps, namely the so called $1$--cone function.
 The location of minima of this function, as well as its values  along a periodic orbit appears to give  valuable information on that orbit.
 \par

 The examples given in this paper
illustrate that a correct characterization of a periodic orbit of a twist map needs the analysis of the residue, the twist number (i.e. the sign of eigenvalues of the  Hessian matrix, $H$, and of its companion matrix $H^-$), as well as the order properties of that orbit.\par
The twist number allows a  better characterization of the periodic orbits given by Aubry--Mather theory:  a minimizing orbit is ordered and has the twist number equal to zero, while a mini-maximizing orbit is also ordered, but its twist number can belong either to the interval $[-1/2, 0)$ or to the interval $(-1, -1/2)$.\par

The study initiated in   this paper allows to investigate  the rotational component of chaos, that along with the Lyapunov component can give a deeper insight into the dynamics of twist area preserving maps.\par

The present study  is  not exhaustive and  implications of the existence of orbits of large twist numbers are to be  studied further. We mention here only few problems that deserve to be investigated:\par
a)  What particular properties of a twist area preserving map  lead to the birth of badly ordered periodic orbits of large absolute twist number.\par
 b) How invariant manifolds of hyperbolic periodic orbits of absolute twist number greater or equal to 1 intersect with the invariant manifolds of  minimizing periodic orbits and how such intersections  influence the transport in the phase space or the dynamics in a region of instability.\par

\section*{Acknowledgements}
This work was  partially supported by Polytechnic University of Timisoara.

\appendix
\section{Translation number of circle homeomorphisms induced by
matrices in $SL(2,\R)$}\label{apendix1}

\subsection{Translation number of an orientation preserving circle homeomorphism}
 Our approach of the twist number of a $q$-periodic orbit of a twist APM exploits some properties of
orientation preserving homeomorphisms (OPHs) of the unit circle, and
of the translation   number of their orbits. For  details concerning
dynamics of OPHs of the circle we refer to \cite{herman, katok,
franks}. \par In the sequel we consider $\mathbb{S}^1$ either as the
unit circle, $C$, in $\R^2$ or as the quotient group $\R/\Z$. The
two definitions are equivalent through the identification:

\beq\label{equivdefscircle}\theta+\Z\in\R/\Z\leftrightarrow
u_\theta=(\cos(2\pi\theta), \sin(2\pi\theta))\in C\subset\R^2\eeq

We denote by $\pi:\R\to\S$, the covering projection,  $\pi(x)=
x\,\,(\mbox{mod}\,\, 1)$. A lift of an OPH $g:\S\to\S$ is a
homeomorphism, $G:\R\to\R$, such that  $\pi\circ G=g\circ \pi$. The
orientation property of $g$ ensures that $G(x+1)=G(x)+1$, $\all\,\,
x\in\R$. Thus the function $\Psi=G-\mbox{id}_\R$ is 1-periodic and
can be written as $\Psi=\psi\circ\pi$, where  $\psi:\S\to\R$. $\psi$
is called
 {\it the displacement map of $g$ associated to $G$} \cite{franks}.  The value of $\psi$ at a point of $\S$ measures
  the amount that point  is displaced by $g$ around the
circle. \par In the sequel we especially use the function
$\Psi=G-\mbox{id}_\R$, and call it  the displacement map of $G$.\par
We  denote by
$m(G)=\min_{x\in\R}\Psi(x)$, $M(G)=\max_{x\in\R}\Psi(x)$.\par Two lifts $G, H$ of $g$ differ by an integer, $H=G+k$,
as well as their displacement maps, $\Psi_H=\Psi_G+k$, $k\in\Z$.\par

The group $\Oph$ of OPHs of the circle has the universal covering group
 \cite{donzhu}:
$$\cOph=\{G:\R\to\R\,\,|\,\, G\,\, \mbox{strictly increasing continuous map,}\,\, G(x+1)=G(x)+1\},$$ i.e. the
elements of this groups are lifts of OPHs of $\S$. We denote by
$\Pi$ the covering projection:
$$\Pi:\widetilde{\mbox{Homeo}^+(\S)}\to\Oph,   \,\, \Pi(G)=g, $$
where $g$ is the unique homeomorphism in $\mbox{Homeo}^+(\S)$ such
that $g\pi=\pi G$. $\Pi$ is a group homomorphism. \par An important
role in the qualitative description of the dynamics of OPHs of the
circle has the Poincar\'{e} translation number. \par The translation number  of a homeomorphism
$G\in\widetilde{\mbox{Homeo}^+(\S)}$ is defined by:
$$\tau(G)=\lim_{n\to\infty}\ds\frac{G^n(x)-x}{n}$$
This number exists   and is independent on $x$. \par

\par For a self-contained presentation  we list some properties of
the homeomorphisms in $\cOph$ and their translation numbers:\par
\begin{prop} \label{proptranslrot} Let $g:\S\to\S$ be an OPH, and
$G$ a lift of $g$.
 Then the following properties hold:\par

 1. The translation number
is invariant under the conjugacy in
$\widetilde{\mbox{Homeo}^+(\S)}$, that is $\tau(H^{-1}\circ G\circ
H^{-1})=\tau(G)$, for every $G,H\in\cOph$.\par
 2. The translation number $\tau: \cOph\to\R$ is not a
homomorphism of groups, that is, in general
$$\tau(G_2 \circ G_1)\neq \tau(G_2)+\tau(G_1),$$ and  its deviation from being a
homomorphism is bounded, existing a positive number $D$ such that:
\beq\label{deffect} |\tau(G_2 \circ G_1)- \tau(G_2)-\tau(G_1)\leq
D,\,\, \all\,\, G_1,G_2 \in \cOph.\eeq The least constant D with
this property is $1$ {\rm\cite{donzhu}}.\par

 3.  $\tau(G^n)=n\tau(G)$,
for any positive integer $n$.\par 4. The translation
number of the inverse of a homeomorphism $G
\in\widetilde{\mbox{Homeo}^+(\S)}$ is $\tau(G^{-1})=-\tau(G)$.\par
5. The difference between the max and min displacement of a lift $G$
is less than 1, $M(G)-m(G)<1$ {\rm\cite{herman}}.\par 6. The translation
number of a homeomorphism, $G$, belongs to the range of the
displacement map, $\tau\in [m(G),M(G)]$.\par
\end{prop}

\subsection{The translation number of the circle homeomorphisms defined  by
matrices in $SL(2,R)$} $\sl2$ is the group of $2\times 2$ real
matrices of determinant equal to $1$. It acts  on the vector space
$\R^2$. This action induces a natural action of the group $\sl2$ on
the unit circle, $\S$, i.e. there is an injective homomorphism
$\phi:\sl2\to \mbox{Homeo}^+(\S)$ between $\sl2$ and the group  of
orientation preserving homeomorphisms of the circle
$C=\{u_\theta=(\cos(2\pi\theta), \sin(2\pi\theta)\in\R^2,
\theta\in\R\}$: \beq\label{delcirclemapA} A\in\sl2\mapsto
g_A\in\mbox{Homeo}^+(\S),\eeq where
$g_A(u_\theta)=\ds\frac{Au_\theta}{\|Au_\theta\|}$, and $\|\cdot\|$
is the Euclidean norm in $\R^2$.\par Therefore the group $\sl2$ is
identified with a subgroup of $\mbox{Homeo}^+(\S)$, and its
universal covering group $ \widetilde{\sl2}$ with a subgroup of
$\widetilde{\mbox{Homeo}^+(\S)}$. We denote by $\Pi:\csl2\to\sl2$
the covering projection.

 Our aim is to
characterize the translation number of the lifts $G\in\csl2$. Its
particular properties in comparison to the translation number of a
general circle homeomorphism derive from the following :
\bpr\label{liftxp12} Any lift $G_A$ of an OPH, $g_A, A\in\sl2$, has
the property \beq\label{GAx12}G_A(x+1/2)=G_A(x)+1/2, \all\,\,
x\in\R\eeq\epr Proof: We denote by $T_\omega$, $\omega\in\R$, the
shift homeomorphism of $\R$, $T_\omega(x)=x+\omega$.  The center of
the group $\sl2$ is $\{\pm I_2\}$.  Because each matrix $A$ commutes
with $-I_2$, the circle map, $g_A$, commutes with the circle
homeomorphism, $h$, defined by:
\beq\label{symmetrygA}h:\theta+\Z\in\R/\Z\mapsto
\theta+1/2+\Z\in\R/\Z\eeq

The lifts of this map are the shifts $T_{(2k+1)/2}$, $k\in\Z$. Thus
a lift $G_A$ commutes with $T_{1/2}$, and this amounts to
$G(x+1/2)=G(x)+1/2$, $\all\,\, x\in\R$.\par\medskip

\par Let  $\mbox{Homeo}^+_h(\S)$ be the subgroup of
orientation preserving homeomorphisms of the circle, that commute
with the map $h$ defined in (\ref{symmetrygA}). Because
$\sl2\subset\mbox{Homeo}^+_h(\S)$, we deduce the   particularities
of the translation number  of the homeomorphisms in
$\mbox{Homeo}^+_h(\S)$ and those properties will be inherited by the
translation numbers of lifts in $\csl2$. The key point in deriving
these particularities is that
 $\mbox{Homeo}^+_h(\S)$ is a rescaled version of
$\widetilde{\mbox{Homeo}^+(\S)}$. Indeed,  let
$\varpi:\widetilde{\mbox{Homeo}^+(\S)}\to\widetilde{\mbox{Homeo}^+(\S)}$
be a map defined by:
\beq\label{maptoantipod}\varpi(H)(x)=\ds\frac{1}{2}H(2x), x\in\R\eeq
$\varpi$ is a monomorphism of groups, and  its image is
$\widetilde{\mbox{Homeo}^+_h(\S)}$.
Hence each $G\in\widetilde{\mbox{Homeo}^+_h(\S)}$ is represented by
a general homeomorphism $H\in\cOph$,  $H=\varpi^{-1}(G)$,
$H(x)=2G(x/2)$.
\begin{prop}\label{translsl2r} Let  $G$ be a homeomorphism in $\widetilde{\mbox{Homeo}^+_h(\S)}$. Then
its translation number is half the translation number of the
corresponding rescaled version $H=\varpi^{-1}(G)$:
$\tau(G)=\tau(H)/2$.\end{prop}

Proof.  From $G(x)=\ds\frac{1}{2}H(2x)$, we get
$G^n(x)=\ds\frac{1}{2}H^n(2x)$, and the translation number is:
$$\tau(G)=\ds\frac{1}{2}\lim_{n\to \infty}\ds\frac{H^n(2x)-2x}{n}=\ds\frac{\tau(H)}{2}$$
\par\medskip

\begin{figure}\label{displacement12}
\centerline{\scalebox{0.75}{\includegraphics{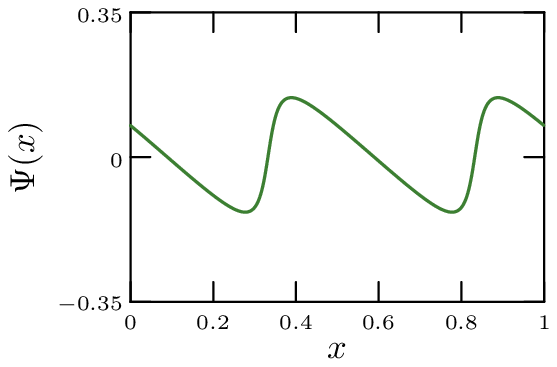}\includegraphics{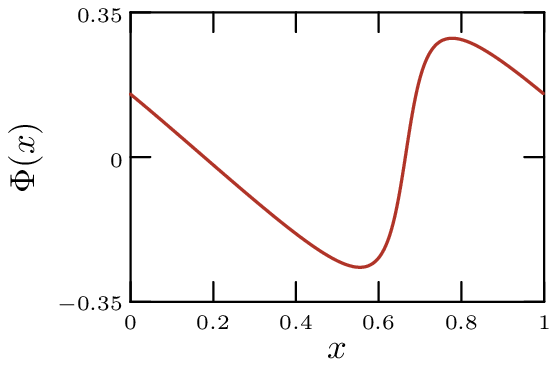}}}
\caption{The graph of displacement maps associated to two
homeomorphisms $G_A, H_A$, $H_A=\varpi^{-1}(G_A)$ .}
\end{figure}

 From the definition of $\varpi$, the Propositions \ref{liftxp12}, \ref{translsl2r}, and the  properties of general circle
homeomorphisms, we have  the following properties of the
 displacement
maps of lifts of circle maps defined by matrices in $\sl2$, relevant
to our approach to twist numbers:
\begin{cor}\label{corrolarytransl}

a) The displacement map of a lift $G_A\in\csl2$ is a periodic
function  of period $1/2$  (Fig.\ref{displacement12}, left).\par
 b) If
$\Psi=G_A-\mbox{id}_\R$ is the displacement map of  a lift
$G_A\in\csl2$, and $\Phi$ is the displacement map of the lift
$H_A=\varpi^{-1}(G_A)$, then $\Psi(x)=\ds\frac{1}{2}\Phi(2x), \all
\,\, x\in\R$ (Fig. \ref{displacement12}). \par
 c) The
difference between the max and min displacement of a homeomorphism
$G_A\in\csl2$ is:
  \beq \label{diffMaxmin} 0\leq
M(G_A)-m(G_A)\leq 1/2\eeq\par
  d) The restriction of
the translation number $\tau$ to the universal covering space of the
group $\sl2$ has the property:\beq\label{quasimsl2} |\tau(G_A\circ
G_B)-\tau(G_A)-\tau(G_B)|\leq 1/2, \quad\all\,\, G_A,
G_B\in\widetilde{\sl2}\eeq
\end{cor}

\section{The structure of the universal covering space of $SL(2,\mathbb{R})$}\label{apendix2}

We define a system of  coordinates on the universal cover group
$\csl2$, that allows us to understand and visualize the topology of this group.

 We recall that the matrices in $\sl2$ are
classified according to their trace.   Namely, a matrix
$A\in\sl2\setminus\{\pm I_2\}$ is called regular hyperbolic, inverse
hyperbolic, elliptic, positive parabolic, respectively negative
parabolic if $\mbox{tr}(A)>2$, $\mbox{tr}(A)<-2$,
$|\mbox{tr}(A)|<2$, $\mbox{tr}(A)=2$, respectively
$\mbox{tr}(A)=-2$.\par

Each matrix $A\in\sl2$ is $\sl2$--conjugated with one of the normal
forms:

\beq \label{normalforms}\barr{l}
R_{2\pi\theta}=\left(\barr{rr}\cos{2\pi\theta}&-\sin{2\pi\theta}\\
\sin{2\pi\theta}&\cos{2\pi\theta}\earr\right),\quad \theta\neq 0, \\
\\H^{r}=\left(\barr{ll} \lam&0\\0&1/\lam\earr\right), \lam\in(0,\infty)\setminus\{1\},
\\  H^i=-H^r\\ P^{1,1}=\left(\barr{rr} 1&1\\0&1\earr\right),\quad
P^{1,-1}=\left(\barr{rr} 1&-1\\0&1\earr\right), \\
P^{-1,1}=\left(\barr{rr} -1&1\\0&-1\earr\right),\quad
P^{-1,-1}=\left(\barr{rr} -1&-1\\0&-1\earr\right)\earr\eeq

It is well known that two conjugated matrices have the same trace.
The conversely is not always true. The matrices within the pairs,
($R_{2\pi\theta}$, $R_{-2\pi\theta}$), ($P^{1,1}, P^{1,-1}$),
($P^{-1,1}, P^{-1,-1}$) for example,  have the same trace, but they
are are not $\sl2$--conjugated.
\par

In the sequel we point out under what conditions a matrix
$A=\left[\barr{ll} a&b\\c&d\earr\right]\in\sl2$ is
$\sl2$--conjugated to one of the above normal forms.

\par
 a) An elliptic matrix, $A\in\sl2$, is $\sl2$--conjugated with a
unique matrix $R_{2\pi\theta}$. More precisely according to whether
$b-c>0$ or $b-c<0$ (conditions that are equivalent respectively to
$b>0$ or $b<0$),
 it is conjugated  to  $R_{2\pi\theta}$,
$\theta\in(-1/2,0)$ or to $R_{2\pi\theta}$, $\theta\in(0,1/2)$ (this property
can be checked  by a straightforward matricial computation). A matrix with
$b=0$ cannot be elliptic.\par Thus $\sl2$ contains two subsets of
elliptic type matrices, one consisting in matrices that are
$\sl2$--conjugated to a clockwise rotation (or equivalently, they
have the entry $b>0$), and another is formed by matrices conjugated
to an anticlockwise rotation (their $b$-entry is negative).

\par b) If $A\in\sl2$ is regular hyperbolic, respectively inverse
hyperbolic, then $A$ is conjugated to an $H^r$, respectively an
$H^{i}$--type matrix. \par

\par c) If $A\in \sl2$ is a positive (negative) parabolic matrix,
then it is conjugated with the first or the second  matrix of the
pairs $(P^{1,1}, P^{1,-1})$, ($P^{-1,1}, P^{-1,-1}$), according to
whether $b-c>0$ or $b-c<0$. \par\medskip

We can conclude that: \par{\it Two hyperbolic matrices are
$\sl2$--conjugated iff they have the same trace. Two elliptic,
 matrices, are $\sl2$--conjugated if and only
if they have the same trace and the same sign of the entries $b$,
and two parabolic matrices are $\sl2$--conjugated iff they have the
same trace  and  the same sign of the differences $b-c$, associated
to each matrix}.\par

Topologically, $SL(2, \R)$ is the interior of a three dimensional
solid torus (\cite{duistermaat}, pag. 14) i.e. it is homeomorphic to
$ D\times \S$, where $D$ is an open disc centered to the origin in
$\R^2$. Hence its universal covering group, $\csl2$, is
homeomorphic(as well as diffeomorphic) to $D\times \R$.   Obviously,
$\csl2$ is also  diffeomorphic to $\R^3$, and a few such
diffeomorphisms are constructed in \cite{burghele}. We define a
diffeomorphism  that fits best to our purpose of studying the twist
number of orbits in area preserving maps. More precisely we choose a
system of coordinates, such that the visual representation  of
$\csl2$ in this system  looks like the Penrose diagram of $\sl2$,
regarded as an anti-de-Sitter space \cite{bengt}.

\par Any matrix $A\in \sl2$ decomposes uniquely as $A=R S$ (polar decomposition), where $R=R_{2\pi\theta}$ is a rotation matrix,
and $S$ is a symmetric positive definite matrix. Let $D$ be the open
disc of radius $1/2$, centered at the origin in $\R^2$, and $(\rho,
\alpha)$ polar coordinates within this disc. We define
 $\chi: D\times\R\to\csl2$, such that
$\chi(0,0,0)=id_\R\in\csl2$, and the polar decomposition of the
matrix $A=\Pi(\chi(\rho, \alpha, \theta))$  be:
\beq\label{polardecomplift} A=R_{2\pi\theta}
R_{2\pi\alpha}\left[\barr{cc} \sqrt{\ds\frac{1-\sin{(\pi
\rho})}{1+\sin({\pi\rho})}}& 0\\0&\sqrt{\ds\frac{1+\sin{(\pi
\rho})}{1-\sin({\pi\rho})}}\earr\right] R_{-2\pi\alpha}\eeq  One can
prove that $\chi$ is a diffeomorphism. Since the trace of the polar
decomposition, $RS$, is
  $\mbox{trace}(R)\mbox{trace}(S)/2$, we get:
   \beq\label{tracecsl2}\mbox{trace}(\Pi(\chi(\theta, \rho, \alpha))=\ds\frac{2\cos{(2\pi\theta})}{\cos(\pi
   \rho)}\eeq

In the sequel the trace of a homemorphism $G\in\csl2$ is meant as
the trace of the corresponding projection, $\Pi(G)$, onto $\sl2$.
According to their traces we also  call the lifts, elliptic,
hyperbolic or parabolic lifts. Denoting by $\mbox{tr}:\csl2\to\R$
the function trace, the connected components of the level sets
$L_t=\mbox{tr}^{-1}\{t\}$, $t\in\R$, are 2-dimensional surfaces and
each such a component represents a conjugacy class in $\csl2$
\cite{burghele}.
\par
\begin{figure}
\centerline{\includegraphics{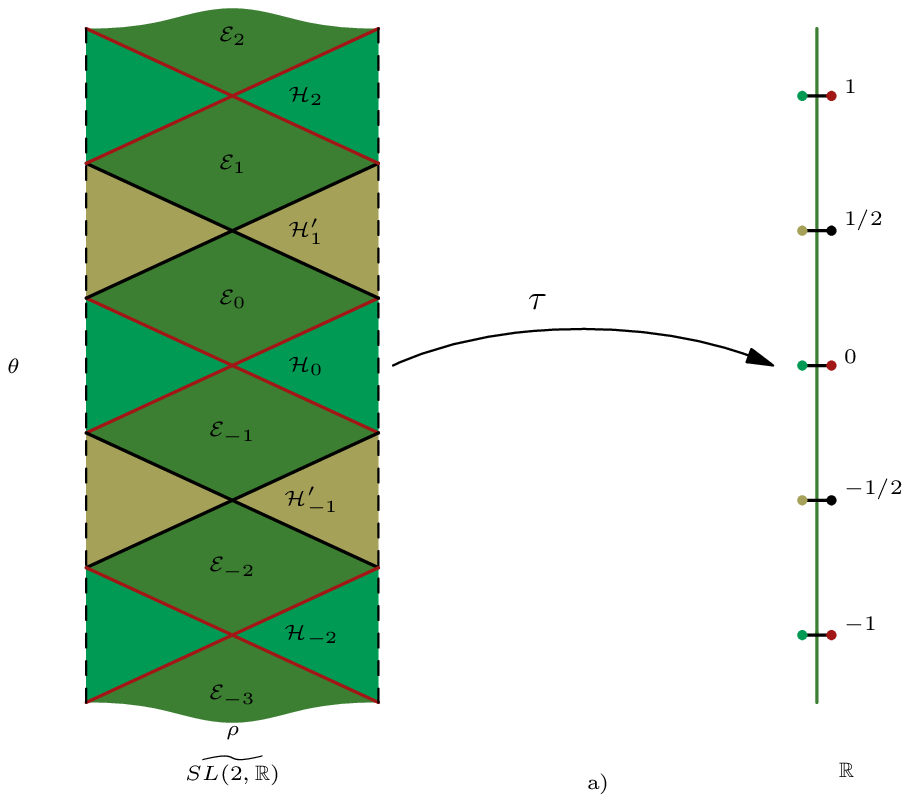}\hskip
2.5cm\includegraphics{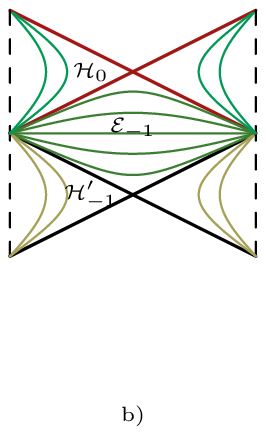}}
\caption{\label{liftsl2r} a) Visualization of the  topological
structure of the group $\csl2$ and  the translation number function
$\tau:\csl2\to\R$; b) The foliation of the regions
$\mathcal{E}_{-1}, \mathcal{H}_0, \mathcal{H}'_{-1}$ by conjugacy
classes.}
\end{figure}

 From (\ref{tracecsl2}) it follows that $\mbox{tr}$ does not depend  on $\alpha$, and thus each  connected
 component of a level set, $L_t$, is a surface of revolution about
 the $\theta$ axis.  This property implies that  the relative position of different
 conjugacy classes
 can be visualized  in any planar section of
the infinite solid cylinder, through a plane containing its symmetry
axis (Fig.\ref{liftsl2r})\par

In order to analyze a  section in the solid cylinder, we  locate the
particular lifts, that project to normal forms in $\sl2$. We denote
by $T_\omega\in\csl2$, the shift homeomorphism, defined by
$T_\omega(x)=x+\omega$. $\Pi(T_\omega)=R_{2\pi\omega}$.\par
 The center of $\csl2$, i.e. the set,
$Z(\csl2)$, of the lifts that commutes with any other lift, consists
in the set of shifts $T_{k/2}$, $k\in\Z$. $\Pi(T_{k/2})=\pm I_2$,
according to as $k$ is  even or odd. $\Pi(\chi(\rho, \alpha, k))$,
$k\in\Z$ is a regular hyperbolic matrix, and $\Pi(\chi(\rho, \alpha,
(2k+1)/2))$, $k\in\Z$, is an inverse hyperbolic matrix.\par

The components of the level set $L_2$, and $L_{-2}$ i.e. the set of
positive, respectively negative  parabolic lifts   are  cones in
$D\times \R$ with the vertex removed.

Analyzing the regions inside the infinite open solid cylinder,
defined by the two families of cones, and the equations of the
conjugacy classes surfaces, $\ds\frac{2\cos{(2\pi\theta})}{\cos(\pi
   \rho)}=t$, corresponding to $t\in (-2,2)$, $t>2$ and $t<-2$ we can conclude that:

\par 1. The parabolic cones have the vertices at a shift $T_{k/2}$,
$k\in\Z$. More precisely, if $k$ is even (odd), then $T_{k/2}$ is
the vertex of a positive (negative) parabolic cone. Each set
$\{T_{k/2}\}$ is a conjugacy class.
\par 2. Inside a region bounded by two patches of nearby cones with
vertices at $T_{k/2}$, $T_{(k+1)/2}$, $k\in\Z$, lie elliptic lifts.
These regions are denoted by $\mathcal{E}_k$. For $k$ even integer,
the matrices $\Pi(\mathcal{E}_k)$ have $b<0$, while for odd $k$ they
have  $b>0$.\par 3. Outside the two sheets of the same cone lie
hyperbolic lifts. For $\theta\in \left(\ds\frac{2k-1}{4},
\ds\frac{2k+1}{4}\right)$, $k$ an even integer, the points $( \rho,
\alpha, \theta)\in D\times \R$, represent regular hyperbolic lifts,
while for $k$ odd integer they represent inverse hyperbolic lifts.
We denoted by $\mathcal{H}_{2j}$, respectively
$\mathcal{H}'_{2j+1}$, $j\in\Z$ the subsets of regular, respectively
inverse hyperbolic lifts.\par Denoting by $\partial (\mathcal{E}_i,
\mathcal{H}_j)$ the boundary between two nearby regions we note
that:

4. The cones $\mathcal{P}_{2k}^{1,1}=\partial
(\mathcal{E}_{2k-1},\mathcal{H}_{2k})\setminus \{T_{k}\}$, $k\in\Z$,
contain lifts of matrices $\sl2$--conjugated to $\left[\barr{ll}
1&1\\0&1\earr\right]$, i.e. matrices with $b-c>0$. The boundary
$\mathcal{P}_{2k}^{1,-1}=\partial(\mathcal{E}_{2k},\mathcal{H}_{2k})\setminus\{T_{k}\},
k\in\Z,$ consists in lifts $G$, whose projections are
$\sl2$--conjugated to $\left[\barr{rr} 1&-1\\0&1\earr\right]$, i.e.
matrices with $b-c<0$.\par

 5. The cones
$\mathcal{P}_{2k-1}^{-1,1}=\partial
(\mathcal{E}_{2k-1},\mathcal{H}'_{2k-1})\setminus \{T_{(2k-1)/2}\}$,
$k\in\Z$, contain lifts of matrices $\sl2$--conjugated to
$\left[\barr{rr} -1&1\\0&-1\earr\right]$, i.e. matrices with
$b-c>0$, while $\mathcal{P}_{2k-1}^{-1,-1}=\partial
(\mathcal{E}_{2k},\mathcal{H}'_{2k+1})\setminus \{T_{(2k-1)/2}\}$,
$k\in\Z$, contains lifts of matrices $\sl2$--conjugated to
$\left[\barr{rr} -1&-1\\0&-1\earr\right]$, i.e. matrices with
$b-c<0$.  \par\medskip

In the study of twist numbers we are interested in the range,
$[m(G), M(G)]$, of the displacement map $\Psi$ of a  lift $G$. In
order to deduce the position of such intervals for different regions
in $\csl2$, we show how are related the intervals $[m(G), M(G)]$,
$[m(H^{-1} G H), M(H^{-1} G H)]$, associated to $\csl2$--conjugated
lifts. \bpr\label{conjugminmax} The minimum, respectively the
maximum displacements, of two  $\csl2$--conjugated homeomorphisms,
 belong to the same demi-unit interval, i.e. an
interval of the form $(k/2, (k+1)/2)$, for some $k\in\Z$: \epr

Proof:

We give the proof only for the minimum displacement of two
$\csl2$--conjugated homeomorphisms. More precisely we show that:

\beq\label{mindisplconj}
\barr{lll}\max\{\ds\frac{n}{2}\in\Z/2\,|\,\ds\frac{n}{2}\leq
m(H^{-1}GH)\}&=&\max\{\ds\frac{n}{2}\in\Z/2\,|\, \ds\frac{n}{2}\leq
m(G)\}\\ & &\\\min\{\ds\frac{n}{2}\in\Z/2\,|\, \ds\frac{n}{2}\geq
m(H^{-1}GH)\}&=&\min\{\ds\frac{n}{2}\in\Z/2\,|\, \ds\frac{n}{2}\geq
m(G)\}\earr\eeq

In order to get  the first relation in (\ref{mindisplconj}) we prove
that for any $G\in\csl2$ we have:

\beq\label{relGxxp}\max\{\ds\frac{n}{2}\in\Z/2\,|\,
\ds\frac{n}{2}\leq x'-x\}=\max\{\ds\frac{n}{2}\in\Z/2\,|\,
\ds\frac{n}{2}\leq G(x')-G(x)\}\eeq

  Let  $x,x'\in\R$ be such that $n/2\leq  x'-x<(n+1)/2$, and distinguish
two cases: i) $x'=x+n/2$,
 and ii) $x'=x+\delta$, $n/2<\delta<(n+1)/2$, $n\in\Z$. In the former case
we have $G(x+n/2)=G(x)+n/2$, i.e. $G(x')-G(x)=n/2$, while in the
latter one, taking into account that $G$ is increasing, we get:
$$
G(x+n/2)< G(x+\delta)<G(x+(n+1)/2), \mbox{i.e.}\quad
n/2<G(x')-G(x)<(n+1)/2$$

In relation  (\ref{relGxxp}) for the lift$H^{-1}$:
$$\max\{\ds\frac{n}{2}\in\Z/2\,|\, \ds\frac{n}{2}\leq
x'-x\}=\max\{\ds\frac{n}{2}\in\Z/2\,|\, \ds\frac{n}{2}\leq
H^{-1}(x')-H^{-1}(x)\},$$
 we take $y=H^{-1}(x)$, and  $x'=G(x)=GH(y)$, and get:

$\max\{\ds\frac{n}{2}\in\Z/2\,|\, \ds\frac{n}{2}\leq G(H(y))-H(y)\}=
\max\{\ds\frac{n}{2}\in\Z/2\,|\, \ds\frac{n}{2}\leq
H^{-1}GH(y)-y\}$, $\all\,\, y\in\R$. Taking the minimum for $y\in\R$
we get the first relation in (\ref{mindisplconj}).\par Similarly one
proves the second one.
\par\medskip

Because each region $\mathcal{E}_k$, $\mathcal{H}_{2k},
\mathcal{H}'_{2k-1}$, $k\in\Z$ is foliated by conjugacy classes, we
have:
 \bpr The displacement intervals for maps $G$ in different regions
of the $\csl2$ are located  as follows:
\par
a) If $G\in \mathcal{E}_k$, $[m(G), M(G)]\subset (k/2,
(k+1)/2)$;\par b) If $G\in \mathcal{H}_{2k}$, then $k\in (m(G),
M(G))$;\par c) If $G\in\mathcal{H}'_{2k+1}$, then
$\ds\frac{2k+1}{2}\in (m(G), M(G))$;\par \epr

Proof. a) A conjugacy class in $\mathcal{E}_k$ is represented by the
shift $T_\omega$, with $\omega\in (k/2, (k+1)/2)$. The displacement
map of $T_\omega$ is the constant function $\Psi(x)=\omega$, and
thus $m(T_\omega)=M(T_{\omega})\in (k/2, (k+1)/2)$. From Proposition
\ref{conjugminmax} we get a).\par b) A conjugacy class in
$\mathcal{H}_{2k}$ is represented by a lift, $G$, of the normal
form, $H^r$ (\ref{normalforms}),
 with $m(G)=k-\delta$, $M(G)=k+\delta$, for some
 $0<\delta<1/4$.\par
 c) A conjugacy class in $\mathcal{H}'_{2k-1}$, $k\in\Z$ is
 represented by a lift $G$ of the normal $H^i$, and  $m(G)=(2k-1)/2-\delta$,
 $M(G)=(2k-1)/2+\delta$, for some $0<\delta<1/4$.\par\medskip
Taking into account the  dynamical properties   of circle maps
\cite{katok}, the properties of the translation numbers stated in
 \ref{apendix1}, and that the translation number of a shift
$T_\omega$ is $\omega$ we get that, if $G\in\csl2$ is a lift of a
circle homeomorphism $g_A$, then the translation number of $G$ is
(see also \cite{mather}, \cite{crovis}):
\begin{enumerate}
\item an integer number, if  $G$ is a regular hyperbolic
or a positive parabolic lift;
\item a demi-integer, i.e. a number in $(2\Z+1)/2$, if $G$ is an inverse
hyperbolic or a negative parabolic lift;
\item a real number different from integers and demi-integers, if
$G$ is an elliptic lift. More precisely $\tau(G)=\theta+k$,
$k\in\Z$, where $\theta\in(-1/2, 0)$ or $(0,1/2)$, depending on the
sign of $b$ in $\Pi(G)$.
\end{enumerate}

In Fig.\ref{liftsl2r} a) we illustrate the values of translation
number function, in each region in the space $\csl2$.

\end{document}